\documentclass[11pt, twoside]{article}

\usepackage{amssymb}
\usepackage{amsmath}

\allowdisplaybreaks

\def\rr{{\mathbb R}}
\def\rn{{{\rr}^n}}
\def\rrm{{{\rr}^m}}

\def\zz{{\mathbb Z}}
\def\nn{{\mathbb N}}

\def\cc{{\mathbb C}}

\def\cx{{\mathcal X}}

\def\cd{{\mathcal D}}

\def\cg{{\mathcal G}}
\def\cl{{\mathcal L}}

\def\ccc{{\mathcal C}}

\def\fz{\infty}
\def\az{\alpha}
\def\supp{{\mathop\mathrm{\,supp\,}}}

\def\loc{{\mathop\mathrm{\,loc\,}}}
\def\lip{{\mathop\mathrm{\,Lip\,}}}
\def\lz{\lambda}
\def\dz{\delta}

\def\ez{\epsilon}
\def\ezl{{\epsilon_1}}

\def\kz{\kappa}
\def\bz{\beta}

\def\gz{{\gamma}}

\def\boz{{\Omega}}

\def\vz{\varphi}
\def\tz{\theta}
\def\sz{\sigma}

\def\wz{\widetilde}

\def\hs{\hspace{0.3cm}}

\def\ls{\lesssim}
\def\gs{\gtrsim}

\def\hl{{H^1_\cl(\rn)}}

\def\bmo{{\mathrm \,BMO\,}}

\def\lip{{\mathrm{\,Lip\,}}}

\def\hl{{\mathop\mathrm {HL}}}

\def\ati{{\mathrm {AOTI}}}
\def\bmo{{\mathop\mathrm{BMO}}}

\def\ocg{{\mathring{\cg}}}

\def\diam{{\mathop\mathrm{\,diam\,}}}

\def\r{\right}
\def\lf{\left}

\newtheorem{thm}{Theorem}[section]
\newtheorem{lem}{Lemma}[section]
\newtheorem{prop}{Proposition}[section]
\newtheorem{rem}{Remark}[section]
\newtheorem{cor}{Corollary}[section]
\newtheorem{defn}{Definition}[section]

\numberwithin{equation}{section}

\begin{document}

\arraycolsep=1pt

\title{\Large\bf Radial maximal function characterizations of
Hardy spaces on RD-spaces and their applications
\footnotetext{\hspace{-0.65cm} Dachun Yang was supported by the National
Natural Science Foundation (Grant No. 10871025) of China.
\noindent$\line(1,0){170}$
\endgraf\noindent D. Yang
\endgraf\noindent School of Mathematical Sciences, Beijing Normal University,
\endgraf\noindent Laboratory of Mathematics and Complex Systems, Ministry of Education,
\endgraf\noindent Beijing 100875, People's Republic of China
\endgraf\noindent e-mail: dcyang@bnu.edu.cn
\medskip
\endgraf\noindent Y. Zhou
\endgraf\noindent School of Mathematical Sciences, Beijing Normal University,
\endgraf\noindent Laboratory of Mathematics and Complex Systems, Ministry of Education,
\endgraf\noindent Beijing 100875, People's Republic of China
\endgraf\noindent e-mail: yuanzhou@mail.bnu.edu.cn
\medskip
\endgraf\noindent Y. Zhou
\endgraf\noindent Department of Mathematics and Statistics, University of Jyv\"askyl\"a,
\endgraf\noindent P. O. Box 35 (MaD), FI-40014, Finland
\endgraf\noindent e-mail: yuzhou@cc.jyu.fi
}}
\author{Dachun Yang and Yuan Zhou}
\date{ }
\maketitle

{\noindent{\bf Abstract}\quad
Let $\cx$ be an RD-space with $\mu(\cx)=\fz$, which means that
$\cx$ is a space of homogeneous type in the sense of Coifman and Weiss
and its measure has the reverse doubling property.
In this paper, we characterize
the atomic Hardy spaces $H^p_{\rm at}(\cx)$
of Coifman and Weiss for $p\in(n/(n+1),1]$ via the radial maximal
function, where $n$ is the ``dimension" of $\cx$, and
the range of index $p$ is the best possible.
This completely answers the question proposed by Ronald R. Coifman
and Guido Weiss in 1977 in this setting,
and improves on a deep result of Uchiyama in 1980
on an Ahlfors $1$-regular space
and a recent result of Loukas Grafakos et al in this setting.
Moreover, we obtain a maximal function theory of localized Hardy spaces
in the sense of Goldberg on RD-spaces by generalizing the
above result to localized Hardy spaces
and establishing the links
between Hardy spaces and localized Hardy spaces.
These results have a wide range of applications. In
particular, we characterize the Hardy spaces
$H^p_{\rm at}(M)$ via the radial maximal function
generated by the heat kernel of the Laplace-Beltrami operator $\Delta$
on complete noncompact connected manifolds $M$ having a doubling property and
supporting a scaled Poincar\'e inequality
for all $p\in(n/(n+\az),1]$, where $\az$ represents the regularity
of the heat kernel. This extends some recent
 results of Russ and Auscher-McIntosh-Russ.}

\medskip



\section{Introduction\label{s1}}

The theory of Hardy spaces on the Euclidean
space $\rn$ plays an important role in various fields of analysis
and partial differential equations; see, for examples,
\cite{sw,fs72,clms,s93,g08}. It is well-known that the following spaces
of homogeneous type of Coifman and Weiss \cite{cw1,cw2} form a natural setting
for the study of function spaces and singular integrals.

\begin{defn}\label{d1.1} \rm
Let $(\cx,\, d)$ be a metric space with a
Borel regular measure $\mu$
such that all balls defined by $d$ have finite and positive measure.
For any $x\in \cx $ and $r>0$, set the ball $B(x,r)\equiv\{y\in \cx :\
d(x,y)<r\}.$ The triple $(\cx,\,d,\,\mu)$ is called a
    space of homogeneous type if there exists
    a constant $C_1\ge 1$ such that for all $x\in \cx$, $\lz\ge1$ and $r>0$,
    \begin{equation}\label{1.1}
    \mu(B(x, \lz r))\le C_1\lz^n\mu(B(x,r)).
    \end{equation}
\end{defn}
Here $n$, if chosen minimal, measures the ``dimension" of the space
$\cx$ in some sense.

Let $\cx$ be a space of homogeneous type as in
Definition \ref{d1.1}. In 1977, Coifman and Weiss \cite{cw2} introduced the
atomic Hardy spaces $H_{\rm at}^p(\cx)$ for $p\in(0, 1]$
and obtained their molecular characterizations. Moreover,
under certain additional geometric condition,
Coifman and Weiss obtained the radial
maximal function characterization of $H_{\rm at}^1(\cx)$.
Then they further asked the following question; see \cite[pp.\,641-642]{cw2}
or \cite[p.\,580]{u80}.

{\it Question 1: Is it possible to characterize $H^p_{\rm at}(\cx)$ for
$p\in (0,\,1]$ in terms of a radial maximal function?}

Recall that an Ahlfors $n$-regular metric measure space
is a space of homogeneous type as in Definition \ref{d1.1}
satisfying that for all $x\in\cx$ and $r\in (0, 2\diam(\cx))$,
$\mu(B(x,r))\sim r^n$, where $\diam(\cx)\equiv\sup_{x,\,y\in\cx}d(x,y)$.
When $\cx$ is an Ahlfors $1$-regular metric measure space,
in 1980, Uchiyama \cite{u80} partially answered this question
by proving the deep result that  for $p\in (p_0,1]$ and functions
in $L^1(\cx)$, the $L^p(\cx)$ quasi-norms of their grand maximal
functions (as in \cite{ms792})
are equivalent to the $L^p(\cx)$ quasi-norms of
their radial maximal functions defined via some kernels in \cite{cw2}.
However, here $p_0>1/2$; see \cite{gly2} for the explicit value of $p_0$.
Also in this setting, when $p\in(1/2, 1]$, Mac\'ias and Segovia \cite{ms792}
characterized  $H_{\rm at}^p(\cx)$ via a grand maximal function.
Observe that when $\cx$ is an Ahlfors $1$-regular metric measure space
and $p\le 1/2$, it is impossible to
characterize $H^p_{\rm at}(\cx)$ via the radial maximal function since
atoms of $H^p_{\rm at}(\cx)$ have only $0$-order vanishing moment.

Recently, the following RD-spaces were introduced in \cite{hmy2},
which are modeled on Euclidean spaces with Muckenhoupt weights,
Ahlfors $n$-regular metric measure spaces (see, for example, \cite{h01}),
Lie groups of polynomial growth
(see, for example, \cite{a, v1, v2})
and Carnot-Carath\'eodory
spaces with doubling measures (see, for example,
\cite{nsw, ns04, s93}).

\begin{defn} \label{d1.2} \rm
The triple $(\cx,\,d,\,\mu)$ is called an RD-space if it is a space
of homogeneous type as in Definition \ref{d1.1} and there
exist constants $C_2>1$ and $C_3>1$ such that for all $x\in \cx$
and $r\in (0, \diam(\cx))$,
\begin{equation}\label{1.2}
\mu(B(x,C_2r))\ge C_3\mu(B(x,r)).\end{equation}
\end{defn}

We point out that the condition \eqref{1.2} can be replaced
by the following geometric condition:  there exists
a constant $a_0>1$ such that for all $x\in\cx$ and
$0<r<\diam(\cx)/a_0$, $B(x,\,a_0r)\setminus B(x,\,r)\neq\emptyset$;
see \cite[Remark 1]{hmy2} and also \cite{yz09} for
some other equivalent characterizations of RD-spaces.
In particular, a connected space of homogeneous type is an RD-space.

Throughout the whole paper, we always assume that
$\cx$ is an RD-space and $\mu(\cx)=\fz$. Recently,
by extending Uchiyama's idea in \cite{u80},
it was proved in \cite{gly2} that there exists a $p_0$ close to $1$
such that when $p\in(p_0,\,1]$, $H^p_{\rm at}(\cx)$ is characterized
by a certain radial maximal function. However, here $p_0>n/(n+1)$.
This partially and affirmatively answers Question $1$
when $\cx$ is an RD-space.
Moreover, a Littlewood-Paley theory for
$H^p_{\rm at}(\cx)$ when $p\in(n/(n+1),1]$ was
established in \cite{hmy1, hmy2} via Calder\'on reproducing formulae.
Also, via inhomogeneous Calder\'on reproducing formulae,
the grand, the nontangential and the dyadic maximal function characterizations
of $H^p_{\rm at}(\cx)$ for $p\in(n/(n+1),1]$ were also obtained in
\cite{gly1}.

On the other hand, let $M$ be a noncompact manifold having
the doubling property and supporting a scaled Poincar\'e inequality and
$\Delta$ the Laplace-Beltrami operator.
There is an increasing interest in the study of
the Hardy spaces and the Riesz transforms on such noncompact manifolds;
see, for example, \cite{acdh,amr,r} and the references therein.
In particular,
let $H^p_\Delta(M)$ be the Hardy space defined
 by the radial maximal function generated
by the heat kernel.
Based on Uchiyama's result \cite{u80},
Russ \cite{r} stated the equivalence of
$H^p_\Delta(M)$ and $H^p_{\rm at}(M)$ for $p\in(p_0,\,1]$ without a proof,
where $p_0$ is not far from $1$;
see also Theorem 8.2 of \cite{amr} for $p=1$.
Then the following question naturally appears.

{\it Question 2: What is the best possible range of
$p\le 1$ which guarantees the equivalence between
$H^p_{\rm at}(M)$ and $H^p_\Delta(M)$.}

This paper is devoted to Questions 1  and 2 above, and extensions
of these results to localized Hardy spaces in the sense of Goldberg (\cite{g}).
{\it Throughout the whole paper, we always assume that $\cx$ is an RD-space
and $\mu(\cx)=\fz$}. (We leave the case when $\mu(\cx)<\fz$
to another forthcoming paper, because of the need to overcome some
additional subtle technical difficulties.)

First, for $p\in(n/(n+1),\,1]$, we characterize
$H^p_{\rm at}(\cx)$ via a certain radial maximal function;
see Theorem \ref{t3.1} and Corollary \ref{c3.1} below.
This completely and affirmatively answers Question 1
of Coifman and Weiss in the case
when $\cx$ is an RD-space and $\mu(\cx)=\fz$.
Observe that when $p\in(0, n/(n+1)]$, the radial maximal function
cannot characterize $H^p_{\rm at }(\cx)$ since its atoms have only
$0$-order vanishing moment. Thus, the range $p\in(n/(n+1),\,1]$
in Theorem \ref{t3.1} and Corollary \ref{c3.1}
is the best possible. Obviously, Theorem \ref{t3.1}
and Corollary \ref{c3.1} generalize the result of
Uchiyama \cite{u80} to RD-spaces and, moreover,
improve the corresponding results in \cite{u80} and
\cite{gly2} by widening the range of the index $p$
to the best possible.
Let $\{S_\ell\}_{\ell\in\zz}$ be
any approximation of the identity (see Definition \ref{d2.1} below).
The proof of Theorem \ref{t3.1} is based on the following two key observations:
 (i)  Any inhomogeneous test function $\phi$ can be decomposed into
$a_{\ell,\phi}S_\ell+ b_{\ell,\phi}\vz$,
where $\vz$ is a homogeneous test function,
$a_{\ell,\phi}$ and $b_{\ell,\phi}$ are constants satisfying certain uniform
estimates in $\ell$ and $\phi$;
 (ii)  Any inhomogeneous distribution uniquely induces a homogeneous
distribution, which enables us to use the homogeneous
Calder\'on reproducing formula
instead of the inhomogeneous one as in \cite{gly1}.
This combined with (i) overcomes the difficulty
caused by the average term appearing in the
inhomogeneous Calder\'on reproducing formula.
Then by some simple calculations, we
control the grand maximal function via the Hardy-Littlewood maximal
function of a certain power of the radial maximal function. This procedure
further implies Theorem \ref{t3.1}, which is much simpler than and
totally different from those used by Fefferman-Stein in \cite{fs72},
Uchiyama in \cite{u80} and Grafakos-Liu-Yang in \cite{gly2}.

Secondly, via a grand maximal function and a variant
of the radial maximal function, we characterize the localized Hardy space
$H^p_{\ell,\,\rm at}(\cx)$ in the sense of Goldberg
with $p\in(n/(n+1),\,1]$ and $\ell\in\zz$;
see Theorems \ref{t3.2} and \ref{t3.4} below.
We point out that the range of the index $p\in(n/(n+1),1]$
here is also the best possible by a reason similar
to the above and that constants appeared in
Theorems  \ref{t3.2} through  \ref{t3.4} are
uniform in $\ell\in\zz$. In fact, let $H^p_\ell(\cx)$ be
the localized Hardy space
defined by the localized grand maximal function.
In Theorem \ref{t3.2} (i), for $p\in(n/(n+1),\,1]$,
we characterize $H^p_\ell(\cx)$
via a variant of the localized radial maximal function,
and for an element in $H^p_\ell(\cx)$, in
Theorem \ref{t3.2} (ii), we further establish the
equivalence between  its  $H^p_\ell(\cx)$-norm  and
the $L^p(\cx)$-norm of its localized radial maximal function.
The proof of Theorem \ref{t3.2} (i)
is also based on the key observation (i) used in the proof of Theorem
\ref{t3.1} and an application of the inhomogeneous Calder\'on reproducing formula.
Due to the inhomogeneity of such Calder\'on reproducing formula,
we only obtain a variant of the localized radial maximal function
characterization for the localized Hardy spaces;
see the extra condition \eqref{3.1} in Theorem \ref{t3.2}.
But this is quite reasonable as the same phenomena happens
in \cite{hmy2} for the localized Littlewood-Paley
characterization of localized Hardy spaces.
The proof of Theorem \ref{t3.2} (ii) requires another
key observation, namely, the size of dyadic cubes
$\{Q_\tau^{\ell,\,\nu}\}_{\tau,\,\nu}$ appearing in the
average term of the inhomogeneous
Calder\'on reproducing formula (see Theorem \ref{t5.2}
below) can be sufficiently small, which allows us
to obtain sufficiently small decay factor determined by
$j$ in the estimate of $J_1$.
This plays a key role in the proof of
Theorem \ref{t3.2} (ii). In Theorem \ref{t3.4},
we establish the equivalence between
the localized Hardy space defined by the grand maximal function and
the one by atoms.
To prove this, we link $H^p_\ell(\cx)$
and $H^p(\cx)$ in Theorem \ref{t3.3}
by using some ideas from \cite{g} and \cite{hmy2}
and the Calder\'on reproducing formula.

Finally, applying Theorems \ref{t3.1} through
\ref{t3.4} to a noncompact manifold
satisfying the doubling property and
supporting a scaled Poincar\'e inequality,
we obtain, in  Proposition \ref{p4.1}, an explicit range  $p\in(n/(n+\az),\,1]$
for Question 2, where $\az\in(0,\,1]$ is the order
of the regularity of the heat kernel.
This range is the best one which can be obtained by the current approach,
but it is not clear if it is optimal.
However, this already extends
Theorem 8.2 of Auscher, McIntosh and Russ \cite{amr}
and the result stated by  Russ in \cite{r}.
We also apply Theorems \ref{t3.1} through
\ref{t3.4} to the Euclidean space $\rrm$ endowed with the measure $w(x)\,dx$,
where  $w\in A_2(\rrm)$ (the class of Muckenhoupt weights),
and to the boundary of an unbounded model
polynomial domain in $\cc^2$ introduced by Nagel and
Stein \cite{ns06} (see also \cite{ns04});
see Propositions \ref{p4.2} and \ref{p4.3} below.
We point out that
Theorems \ref{t3.1} through \ref{t3.4}
are also valid for Lie groups of polynomial growth.
On the other hand, in this setting, Saloff-Coste has already obtained
certain grand and radial maximal functions
and atomic characterizations
of Hardy spaces;
see \cite{s86,s87,s89,s90} for more details.

The paper is organized as follows.  We recall, in Section 2,  some
notation and definitions;
state, in Section 3, the main results of
this paper, Theorems \ref{t3.1} through \ref{t3.4};
give, in  Section 4, some applications;
and finally, in Section 5, prove Theorems
\ref{t3.1} through \ref{t3.4} by employing
Calder\'on reproducing formulae (see Theorems \ref{5.1} and \ref{t5.2})
 established in \cite{hmy2, gly2}.

We finally make some conventions. Throughout this paper, we always
use $C$ to denote a positive constant which is independent of the
main parameters involved but whose value may differ from line to
line. Constants with subscripts do not change throughout the whole
paper. Set $ a\wedge b\equiv\min\{a,\,b\}$ for $a,\,b\in\rr$.

\section{Preliminaries\label{s2}}

The following notion of
approximations of the identity
on RD-spaces were first introduced in \cite{hmy2}.
In what follows, we set $V_r(x)\equiv\mu(B(x,\,r))$ and
$V(x,\,y)\equiv\mu(B(x,\,d(x,\,y)))$ for $x,\,y\in\cx$ and $r\in(0,\,\fz)$.

\begin{defn} \label{d2.1}\rm
Let $\ez_1\in(0,\,1]$, $\ez_2>0$ and $\ez_3>0$. A sequence
$\{S_k\}_{k\in\zz}$ of bounded linear integral operators on
$L^2(\cx)$ is called an approximation of the identity of order
$(\ez_1,\,\ez_2,\,\ez_3)$ (for short,
$(\ez_1,\,\ez_2,\,\ez_3)$-$\ati$), if there exists a positive
constant $C_4$ such that for all $k\in\zz$ and $x,\, x',\,
y,\,y'\in\cx$, $S_k(x,y)$, the integral kernel of $S_k$, is a
measurable function from $\cx\times\cx$ into $\cc$ satisfying
\begin{enumerate}
\vspace{-0.3cm}
\item[(i)] $|S_k(x,y)|\le C_4\frac 1{V_{2^{-k}}(x)+V(x,\,y)}
[\frac{2^{-k}}{2^{-k}+d(x,\,y)}]^{\ez_2};$ \vspace{-0.3cm}
\item[(ii)] $|S_k(x,y)-S_k(x',y)|
\le C_4 [\frac{d(x,\,x')}{2^{-k}+d(x,\,y)}]^\ezl
\frac
1{V_{2^{-k}}(x)+V(x,\,y)}[\frac{2^{-k}}{2^{-k}+d(x,\,y)}]^{\ez_2}$
for $d(x,x')\le[2^{-k}+d(x,\,y)]/2;$ \vspace{-0.3cm}
\item[(iii)] Property (ii) also holds with $x$ and $y$
interchanged; \vspace{-0.3cm}
\item[(iv)] $|[S_k(x,y)-S_k(x,y')]-[S_k(x',y)-S_k(x',y')]|\le C_4
[\frac{d(x,\,x')}{2^{-k}+d(x,\,y)}]^\ezl
[\frac{d(y,\,y')}{2^{-k}+d(x,\,y)}]^\ezl$
\newline
$\times\frac 1{V_{2^{-k}}(x)+V(x,\,y)}
[\frac{2^{-k}}{2^{-k}+d(x,\,y)}]^{\ez_3}$ for
$d(x,x')\le[2^{-k}+d(x,\,y)]/3$ and $d(y,y')\le[2^{-k}+d(x,\,y)]/3;$
\vspace{-0.3cm}
\item[(v)] $\int_\cx S_k(x,z)\,d\mu(z)=1=\int_\cx S_k(z,y)\,d\mu(z)$
for all $x, \, y\in\cx$.
\end{enumerate}
\end{defn}

\begin{rem}\label{r2.2}\rm
(i) In \cite{hmy2}, for any $N>0$, it was proved that there exists
a $(1,\,N,\,N)$-$\ati$ $\{S_k\}_{k\in\zz}$ with bounded support in the
sense that $S_k(x,y)=0$ when $d(x,\,y)> C2^{-k}$, where $C$ is a
fixed positive constant independent of $k$. In this case,
$\{S_k\}_{k\in\zz}$ is called a $1$-$\ati$ with bounded support;
see \cite{hmy2}.

(ii) If a sequence $\{\wz S_t\}_{t>0}$ of bounded linear integral
operators on $L^2(\cx)$ satisfies (i) through (v) of Definition
\ref{d2.1} with  $2^{-k}$ replaced by $t$, then we call $\{\wz
S_t\}_{t>0}$ a continuous approximation of the identity of
order $(\ez_1,\,\ez_2,\,\ez_3)$ (for short, continuous
$(\ez_1,\,\ez_2,\,\ez_3)$-$\ati$). For example, if
$\{S_k\}_{k\in\zz}$ is an $(\ez_1,\,\ez_2,\,\ez_3)$-$\ati$ and if we
set $\wz S_t(x,\,y)\equiv S_k(x,\,y)$ for $t\in(2^{-k-1},\,2^{-k}]$
with $k\in\zz$, then  $\{\wz S_t\}_{t>0}$ is a continuous
$(\ez_1,\,\ez_2,\,\ez_3)$-$\ati$.

(iii) If $S_k$ (resp. $\wz S_t$) satisfies (i), (ii), (iii) and (v)
of Definition \ref{d2.1}, then $S_kS_k $ (resp. $\wz S_t\wz S_t$)
satisfies the conditions (i) through (v) of Definition \ref{d2.1};
see \cite{hmy2}.

(iv) For any RD-space $(\cx,\,d,\,\mu)$,
if we relax $d$ to be a  quasi-metric,
then there exist constants $\tz\in(0,\,1)$ and $C>0$ and
quasi-metric $\wz d$ which is equivalent to $d$
such that $|\wz d(x,\,y)-\wz d(z,\,y)|\le
 C [\wz d(z,x)]^\tz[\wz d(x,\,y)+\wz d(z,\,y)]$
 for all $x,\,y,\,z\in\cx$; see \cite{ms791}.
By this and an argument  similar to that used in \cite{hmy2},
we know the existence of the approximation of the identity
  $(\tz,\,N,\,N)$-$\ati$ with bounded support, where $N>0$.
\end{rem}

The following spaces of test functions
play an important role in the theory of
function spaces on spaces of homogeneous type; see \cite{hmy1, hmy2}.

\begin{defn} \label{d2.2}\rm
Let $x\in\cx$,  $r>0$, $\bz\in(0,\,1]$ and $\gz>0$. A function $f$
on $\cx $ is said to belong to the space of test functions,
$\cg(x,\,r,\,\bz,\,\gz)$, if there exists a nonnegative constant $C$ such that
\begin{enumerate}
    \vspace{-0.3cm}
    \item[(i)] $|f(y)|\le C\frac 1{V_r(x)+V(x,\,y)}
    \lf(\frac r{r+d(x,\,y)}\r)^\gz$ for all $y\in\cx$;
    \vspace{-0.3cm}
    \item[(ii)] $|f(z)-f(y)|\le C\lf(\frac {d(y,\,z)}{r+d(x,\,y)}\r)^\bz
    \frac 1{V_r(x)+V(x,\,y)}\lf(\frac r{r+d(x,\,y)}\r)^\gz$
    for all $y,\, z\in \cx $ satisfying that $d(y,\,z)\le [r+d(x,\,y)]/2$.
    \vspace{-0.3cm}
\end{enumerate}
Moreover, for any $f\in \cg(x,\,r,\,\bz,\,\gz)$, we define its norm by
$$\|f\|_{\cg(x,\,r,\,\bz,\,\gz)}\equiv\inf\lf\{C:\, (i)\ \mathrm{and}\ (ii)\
\mathrm{hold}\r\}.$$

The space $\ocg(x,\,r,\,\bz,\,\gz)$ is defined to be the set of all functions
$f\in\cg(x,\,r,\,\bz,\,\gz)$ satisfying that $\int_\cx f(y)\,d\mu(y)=0.$
Moreover, we endow the space $\ocg(x,\,r,\,\bz,\,\gz)$ with the same norm as
the space $\cg(x,\,r,\,\bz,\,\gz)$.
\end{defn}

It is easy to see that $\cg(x,\,r,\,\bz,\,\gz)$ is a Banach space.
Let $\ez\in(0,\,1]$ and $\bz,\,\gz\in(0,\,\ez]$. For applications,
we further define the space $\cg_0^\ez(x,\,r,\,\bz,\,\gz)$ to be the
completion of the set $\cg(x,\,r,\,\ez,\,\ez)$ in
$\cg(x,\,r,\,\bz,\,\gz)$. For $f\in\cg_0^\ez(x,\,r,\,\bz,\,\gz)$,
define
$\|f\|_{\cg_0^\ez(x,\,r,\,\bz,\,\gz)}\equiv\|f\|_{\cg(x,\,r,\,\bz,\,\gz)}$.
Let $\lf(\cg_0^\ez(x,\,r,\,\bz,\,\gz)\r)'$ be the set of all
continuous linear functionals on $\cg_0^\ez(x,\,r,\,\bz,\,\gz)$, and
as usual, endow $\lf(\cg_0^\ez(x,\,r,\,\bz,\,\gz)\r)'$ with the weak$^\ast$ topology.
Throughout the whole paper, we fix $x_1\in\cx$ and
write $\cg(\bz,\,\gz)=\cg(x_1,\,1,\,\bz,\,\gz)$, and
$\lf(\cg_0^\ez(\bz,\,\gz)\r)'\equiv\lf(\cg_0^\ez(x_1,\,1,\,\bz,\,\gz)\r)'$.
Observe that for any $x\in\cx$ and $r>0$,
$\cg_0^\ez(x,\,r,\,\bz,\,\gz)=\cg_0^\ez(\bz,\,\gz)$ with equivalent norms.

Similarly, define the space $\ocg_0^\ez(x,\,r,\,\bz,\,\gz)$ to be
the completion of the set $\ocg(x,\,r,\,\ez,\,\ez)$ in
$\ocg(x,\,r,\,\bz,\,\gz)$. For
$f\in\ocg_0^\ez(x,\,r,\,\bz,\,\gz)$,  define
$\|f\|_{\ocg_0^\ez(x,\,r,\,\bz,\,\gz)}\equiv\|f\|_{\ocg(x,\,r,\,\bz,\,\gz)}$.
Denote by $(\ocg_0^\ez(x,\,r,\,\bz,\,\gz))'$ the set of all
continuous linear functionals from
 $\ocg_0^\ez(x,\,r,\,\bz,\,\gz)$ to $\cc$,
 and endow $(\ocg_0^\ez(x,\,r,\,\bz,\,\gz))'$ with the
weak$^\ast$ topology. Write
$\ocg_0^\ez(\bz,\,\gz)\equiv\ocg(x_1,\,1,\,\bz,\,\gz)$. For any
$x\in\cx$ and $r>0$, we also have
$\ocg_0^\ez(x,\,r,\,\bz,\,\gz)=\ocg_0^\ez(\bz,\,\gz)$ with equivalent
norms.

Now we recall the following maximal functions.
\begin{defn}\label{d2.3}\rm
(i) Let $\ez\in(0,\,1]$, $\bz,\,\gz\in(0,\,\ez)$ and $\ell\in\zz$.
For any $f\in
(\cg_0^\ez(\bz,\,\gz))'$, the grand maximal function
$G^{(\ez,\,\bz,\,\gz)}(f)$ is defined by setting, for all $x\in\cx$,
$$
 G^{(\ez,\,\bz,\,\gz)}(f)(x)\equiv\sup\lf\{\langle f,\,\varphi\rangle:\,
\varphi\in\cg^\ez_0(\bz,\,\gz),\
\|\vz\|_{\cg(x,\,r,\,\bz,\,\gz)}\le1 \ \mbox{for some}\ r>0\r\},
$$
and the localized grand maximal function
$G^{(\ez,\,\bz,\,\gz)}_\ell(f)$ by setting, for all $x\in\cx$,
$$
 G^{(\ez,\,\bz,\,\gz)}_\ell(f)(x)\equiv\sup\lf\{\langle f,\,\varphi\rangle:\,
\varphi\in\cg^\ez_0(\bz,\,\gz),\
\|\vz\|_{\cg(x,\,r,\,\bz,\,\gz)}\le1\ \mbox{for some}\ r\in(0,\,2^{-\ell}]\r\}.
$$

(ii) Let $\ez_1\in(0,\,1]$,\,$\ez_2,\,\ez_3>0$,
$\ez\in(0,\,\ez_1\wedge\ez_2)$ and $\{S_k\}_{k\in\zz}$ be an
$(\ez_1,\,\ez_2,\,\ez_3)$-$\ati$. Let $\ell \in\zz$.
For any $\bz,\,\gz\in(0,\,\ez)$ and
$f\in (\cg_0^\ez(\bz,\,\gz))'$, the radial maximal function $
S^+(f)$ is defined by setting, for all $x\in\cx$,
$$
 S^+(f)(x)\equiv\sup_{k\in\zz} |S_k(f)(x)|,
$$
and the localized radial maximal function
$S^+_\ell(f)$ by setting, for all $x\in\cx$,
$$
  S^+_\ell(f)(x)\equiv\sup_{k\ge\ell} |S_k(f)(x)|.
$$
\end{defn}

When there exists no ambiguity, we write
$G^{(\ez,\,\bz,\,\gz)}(f)$ and
$G^{(\ez,\,\bz,\,\gz)}_\ell(f)$ simply by $G(f)$ and $G_\ell(f)$, respectively.
It is easy to see that for all $x\in\cx$, $S^+(f)(x)\le G(f)(x)$,
$S^+_{\ell}(f)(x)\le G_{\ell}(f)(x)$ and
for any $\ell\ge k$,
$G_\ell(f)(x)\le G_k(f)(x)\le C G_\ell(f)(x)$,
where $C$ is a positive constant depending on $k$
and $\ell$, but not on $f$ and $x$.

\begin{defn}\label{d2.4}\rm
Let $p\in(n/(n+1),\,1]$ and $n(1/p-1)<\bz,\,\gz<\ez<1$.

(i) The Hardy space
$H^p(\cx)$ is defined by
$$H^p(\cx)\equiv\lf\{f\in(\cg_0^\ez(\bz,\,\gz))': \
\|f\|_{H^p(\cx)}\equiv\|G(f)\|_{L^p(\cx)}<\fz\r\}.$$

(ii) Let $\ell\in\zz$. The localized Hardy space
$H^p_\ell(\cx)$ is defined by
$$H^p_\ell(\cx)\equiv\lf\{f\in(\cg_0^\ez(\bz,\,\gz))': \
\|f\|_{H^p_\ell(\cx)}\equiv\|G_\ell(f)\|_{L^p(\cx)}<\fz\r\}.$$

\end{defn}

It was proved in \cite{gly1} that if $p\in(n/(n+1),\,1]$, then
the definition of $H^p(\cx)$ is independent of
the choices of  $ \ez\in(n(1/p-1),\,1)$ and $\bz,\,\gz\in(n(1/p-1),\,\ez)$.
This also holds for $H^p_\ell(\cx)$ by a similar argument. Here we omit the details.

Now we recall the notion of the atomic Hardy space of Coifman and Weiss \cite{cw2}.

\begin{defn}\label{d2.5}\rm
Let $p\in(0,\,1]$, $q\in[1,\,\fz]\cap (p,\,\fz]$ and $\ell\in\zz$.

\noindent (i) A measurable function $a$ is called a $(p,\,q)$-atom
associated to the ball $B(x,\,r)$ if

(A1) $\supp a\subset B(x,\,r)$ for certain $x\in\cx$ and $r>0$,

(A2) $\|a\|_{L^q(\cx)}\le [\mu(B(x,\,r))]^{1/q-1/p}$,

(A3) $\int_\cx a(x)\,d\mu(x)=0$.

\noindent (ii) A measurable function $a$ is called a
$(p,\,q)_\ell$-atom
associated to the ball $B(x,\,r)$ if $r\le 2^{-\ell}$ and
$a$ satisfies (A1) and (A2),
and when $r<2^{-\ell}$, $a$ also satisfies (A3).
\end{defn}

\begin{defn}\label{d2.6}\rm
Let $p\in(0,\,1]$.

(i) The space $\lip(1/p-1,\,\cx)$  is defined to be the collection of
all functions $f$  satisfying
$$\|f\|_{\lip(1/p-1,\,\cx)}\equiv\sup_{x,\,y\in\cx,\,B\ni x,\,y}
 [\mu(B)]^{1-1/p} |f(x)-f(y)| <\fz,$$
where the supremum is taken over all $x,\,y\in\cx$
and all balls containing $x$ and $y$.

(ii) The space $\lip_\ell(1/p-1,\,\cx)$ is defined to be the collection of
all functions $f$  satisfying
$$\|f\|_{\lip_\ell(1/p-1,\,\cx)}\equiv\sup_{x,\,y\in\cx,\,B\in I_\ell(x,\,y)}
\frac{|f(x)-f(y)|}{[\mu(B)]^{1/p-1}}+
\sup_{B,\, r_B>2^{-\ell}}\frac1{\mu(B)}\int_{B}|f(z)|\,d\mu(z)<\fz,$$
where $I_\ell(x,\,y)$ denotes all balls containing $x$ and $y$
with radius no more than $2^{-\ell}$,
the second supremum is taken over all balls with radius more than $2^{-\ell}$.

\end{defn}

\begin{defn}\label{d2.7}\rm
Let $p\in(0,\,1]$ and $q\in[1,\,\fz]\cap(p,\,\fz]$.

(i) The space $H^{p,\,q}(\cx)$ is defined to be the
set of all $f=\sum_{j\in\nn}\lz_ja_j$ in $(\lip(1/p-1,\,\cx))'$ when $p<1$ and
in $L^1(\cx)$ when $p=1$,
where $\{a_j\}_{j\in\nn}$ are $(p,\,q)$-atoms
and $\{\lz_j\}_{j\in\nn}\subset\cc$
such that $\sum_{j\in\nn}|\lz_j|^p<\fz$.
For any $f\in H^{p,\,q}(\cx)$,
define $\|f\|_{H^{p,\,q}(\cx)}
\equiv\inf\{(\sum_{j\in\nn}|\lz_j|^p)^{1/p}\}$, where
the infimum is taken over all the above decompositions of $f$.

(ii) The space $H^{p,\,q}_\ell(\cx)$ is defined
as in (i) with $(p,\,q)$-atoms replaced by $(p,\,q)_\ell$-atoms
and  $(\lip(1/p-1,\,\cx))'$ replaced by  $(\lip_\ell(1/p-1,\,\cx))'$.
\end{defn}

Since $H^{p,\,q}(\cx)= H^{p,\,\fz}(\cx)$ and
$H_\ell^{p,\,q}(\cx)= H_\ell^{p,\,\fz}(\cx)$
(see \cite{cw2} and also \cite{yyz08}),
we always write $H^{p,\,q}(\cx)$ and $H_\ell^{p,\,q}(\cx)$
as $H^{p}_{\rm at}(\cx)$ and $H^{p}_{\ell,\,\rm at}(\cx)$,
respectively.
Moreover, the dual spaces of
$H^p_{\rm at}(\cx)$ and  $H^p_{\ell,\,\rm at}(\cx)$
are, respectively, $\lip(1/p-1,\,\cx)$ and $\lip_\ell(1/p-1,\,\cx)$ when $p<1$,
and $\bmo(\cx)$ and $\bmo_\ell(\cx)$ when $p=1$;
see \cite{cw2} and also \cite{yyz08} for the details.

\section{Main results}\label{s3}

The first result is on the
characterization of the radial maximal function of the
Hardy space $H^p(\cx)$.

\begin{thm}\label{t3.1}
Let $\ez_1\in(0,\,1]$,\,$\ez_2,\,\ez_3>0$, $\ez\in(0,\,\ez_1\wedge\ez_2)$
and $\{S_k\}_{k\in\zz}$ be an
$(\ez_1,\,\ez_2,\,\ez_3)$-$\ati$.
Let $p\in(n/(n+\ez),\,1]$ and $\bz,\,\gz\in(n(1/p-1),\,\ez)$.
Then for any $f\in(\cg^\ez_0(\bz,\,\gz))'$, $f\in H^p(\cx)$ if and only if
$\|S^+(f)\|_{L^p(\cx)}<\fz$;
moreover, for all $f\in H^p(\cx)$,
$ \|f\|_{H^p(\cx)}\sim\|S^+(f)\|_{L^p(\cx)}.$
\end{thm}

The proof of Theorem \ref{t3.1} is given in Section \ref{s5}.
 Moreover, by Theorem 4.16 in \cite{gly1} and Theorem \ref{t3.1},
 we have the following result.

\begin{cor} \label{c3.1}
Let $\ez_1\in(0,\,1]$,\,$\ez_2,\,\ez_3>0$, $\ez\in(0,\,\ez_1\wedge\ez_2)$ and
$\{S_k\}_{k\in\zz}$ be an $(\ez_1,\,\ez_2,\,\ez_3)$-$\ati$.
Let $p\in(n/(n+\ez),\,1]$ and $\bz,\,\gz\in(n(1/p-1),\,\ez)$.
Then $f\in H^p(\cx)$ if and only if $f\in H^{p}_{\rm at}(\cx)$ or if and only if
$f\in (\cg^\ez_0(\bz,\,\gz))'$ and
$ S^+(f)\in L^p(\cx)$;
moreover, for all $f\in H^{p}_{\rm at}(\cx)$,
$ \|S^+(f)\|_{L^p(\cx)}\sim\|f\|_{H^{p}(\cx)}\sim\|f\|_{H_{\rm at}^p(\cx)}.$
\end{cor}

\begin{rem} \rm\label{r3.1}

(i) If $\{S_k\}_{k\in\zz}$ is replaced by a continuous $(\ez_1,\,
\ez_2,\,\ez_3)$-$\ati$ as in Remark \ref{r2.2} (ii),
then Theorem \ref{t3.1} and Corollary \ref{c3.1} still hold for all $p\in(n/(n+\ez_1),\,1]$.

(ii) If we relax $d$ to be a quasi-metric,
then Theorem \ref{t3.1} and Corollary \ref{c3.1} still
hold by replacing $p\in(n/(n+1),1]$
with $p\in(n/(n+\tz),\,1]$
 for any $( \tz,\,
1,\,1)$-$\ati$  $\{S_k\}_{k\in\zz}$,
where $\tz\in(0,\,1)$ is the same as in Remark \ref{r2.2} (iv).

(iii) Corollary \ref{c3.1} tells us that for $p\in(n/(n+1),\,1]$,
$H^p_{\rm at}(\cx)$
is characterized by the radial maximal function
in Definition \ref{d2.3} (ii), which completely answers Question 1
when $\mu(\cx)=\fz$
asked by Coifman and Weiss. We also remark that
Theorem \ref{t3.1} and Corollary \ref{c3.1}
improve the deep result of Uchiyama \cite{u80} and that of
\cite{gly2} to the best range $p\in (n/(n+1),\,1]$.

(iv) Notice that atoms of $H^p_{\rm at}(\cx)$ have only $0$-order vanishing moment
for $p\in (0,\,n/(n+1)]$ here.
It is easy to see that the Poisson kernel is just a $(1,\,1,\,1)$-$\ati$
as in Definition \ref{d2.1},
and the radial Poisson maximal function characterizes
certain atomic Hardy spaces on $\rn$, which
asks the vanishing moment of atoms no less than
$1$-order, when $p\in (0,\,n/(n+1)]$;
see, for example, \cite[pp.\,91,\,107,\,133]{s93}.
This atomic Hardy space is essentially different
from $H^p_{\rm at}(\cx)$ considered here.
So it is impossible to use the Poisson maximal function to
characterize  $H^p_{\rm at}(\cx)$ when $p\in (0,\,n/(n+1)]$
here. In this sense, the range
$p\in(n/(n+1),\,1]$ is the best possible for which
$H^{p}_{\rm at}(\cx)$ can be characterized by the radial maximal function.
\end{rem}

Now we turn to the localized Hardy spaces.

\begin{thm}\label{t3.2}
Let $\ez_1\in(0,\,1]$,\,$\ez_2,\,\ez_3>0$,
$\ez\in(0,\,\ez_1\wedge\ez_2)$ and $\{S_k\}_{k\in\zz}$ be an
$(\ez_1,\,\ez_2,\,\ez_3)$-$\ati$.
Let $\ell\in\zz$, $p\in(n/(n+\ez),\,1]$ and $\bz,\,\gz\in(n(1/p-1),\,\ez)$.

(i) Then
$f\in H^p_\ell(\cx)$ if and only if
$f\in(\cg^\ez_0(\bz,\,\gz))'$,
$S^{+}_{k+1}(f)\in L^p(\cx)$, and for any $k\in\zz$ and $a>0$,
 \begin{eqnarray}\label{3.1}
S^{(a)}_k(f)(x)\equiv
\frac{1}{V_{a2^{-k}}(x)}\int_{B(x,\,a2^{-k})}|S_{k}(f)(y)|\,d\mu(y)\in  L^p(\cx),
\end{eqnarray}
where $x\in\cx$.
Moreover, for any $a>0$, there exists a positive
constant  $C$  depending on $a$   such that for all
 $\ell\in\zz$ and $f\in(\cg^\ez_0(\bz,\,\gz))'$,
$$C^{-1} \|f\|_{H^p_\ell(\cx)}
 \le \|S^{+}_{\ell+1}(f)\|_{L^p(\cx)}+ \|S^{(a)}_\ell(f)\|_{L^p(\cx)}
\le C\|f\|_{H^p_\ell(\cx)}.$$

(ii) If $f\in(\cg^\ez_0(\bz,\,\gz))'$ and \eqref{3.1} holds for any $k\in\zz$, then
$f\in H^p_\ell(\cx)$ if and only if
$S^{+}_{\ell}(f)\in L^p(\cx)$; moreover,
there exists a positive constant  $C$, independent of
 $k$ and $\ell$,  such that for all $f\in H^p_\ell(\cx)$,
$C^{-1} \|f\|_{H^p_\ell(\cx)}
 \le \|S^{+}_\ell(f)\|_{L^p(\cx)}
\le C\|f\|_{H^p_\ell(\cx)}.$
\end{thm}

The basic idea of the
proof of Theorem \ref{t3.2} (i) is similar to that used
in the proof of Theorem \ref{t3.1}.
To prove Theorem \ref{t3.2} (ii),
observe that  \eqref{3.1} and $S^+_k(f)\in L^p(\cx)$ imply that
$G_\ell(f)\in L^p(\cx)$
by Theorem \ref{t3.2} (i).
Based on the observation that
the constants appearing in Calder\'on reproducing formulae are uniform in $j$,
where $j$ measures the size of dyadic cubes $\{Q^{k,\,\nu}_{\tau }\}_{\tau,\,\nu}$
appearing in the average term of
the inhomogeneous Calder\'on reproducing formula (see Theorem \ref{t5.2} below).
Then we prove that for certain
$r\in(n/(n+\ez_1),\,p)$, and all $j$ large enough and $x\in\cx$,
$$G_\ell(f)(x)\le C  2^{jn(1/r-1)}[\hl([S^+_\ell(f)]^r)(x)]^{1/r}
+C2^{j(\ez_1+n-n/r)}[\hl([G_\ell(f)]^r)(x)]^{1/r} ,$$
where  $C$ is a positive constant independent of $j$, $f$ and $x$
and  $\hl$ denotes the Hardy-Littlewood maximal function.
Taking $j$ such that
$C2^{j(\ez_1+n-n/r)}\le 1/2$,
we then obtain Theorem \ref{t3.2} (ii).

By first establishing a connection between $H^p(\cx)$ and $H^p_{\ell}(\cx)$,
  we then obtain the
equivalence between $H^p_\ell(\cx)$ and $H^p_{\ell,\,\rm at}(\cx)$
via Corollary \ref{c3.1}.
\begin{thm}\label{t3.3}
Let $\ez_1\in(0,\,1]$,\,$\ez_2,\,\ez_3>0$, $\ez\in(0,\,\ez_1\wedge\ez_2)$ and
$\{S_k\}_{k\in\zz}$ be an $(\ez_1,\,\ez_2,\,\ez_3)$-$\ati$.
Let $p\in(n/(n+\ez),\,1]$ and $\bz,\,\gz\in(n(1/p-1),\,\ez)$.
Then there exists a positive constant $C$ such that
for all $\ell \in\zz$ and $f \in (\cg^\ez_0(\bz,\,\gz))'$,
$\|f- S_\ell(f)\|_{H^p(\cx)}\le C\|G_\ell(f)\|_{L^p(\cx)}.$
\end{thm}

\begin{thm} \label{t3.4}
Let $\ez\in(0,\,1]$, $p\in(n/(n+\ez),\,1]$ and $\bz,\,\gz\in(n(1/p-1),\,\ez)$.
Then for each $\ell\in\zz$,
$f\in H_{\ell,\,\rm at}^{p}(\cx)$ if and only if
$f\in (\cg^\ez_0(\bz,\,\gz))'$  and
$\|G_\ell(f)\|_{L^p(\cx)}<\fz$;
 moreover, there exists a positive constant $C$,
independent of $\ell$,
such that for all $f\in H^p_{\rm at}(\cx)$,
$C^{-1}\|f\|_{H_{\ell,\,\rm at}^{p}(\cx)}\le
\|G_\ell(f)\|_{L^p(\cx)}\le C\|f\|_{H_{\ell,\,\rm at}^{p}(\cx)}.
$
\end{thm}

\begin{rem}\label{r3.2}\rm
(i) If $p=1$, then Theorem \ref{t3.2} (i) gives
 the radial maximal function characterization
of $H^p_\ell(\cx)$ since \eqref{3.1} is just $S_\ell(f)\in L^1(\cx)$.
This result when $p=1$ was also obtained in \cite{yz08}.
When $p<1$, it is still unclear if one can remove the extra assumption
\eqref{3.1}. But this is quite reasonable since we use the
inhomogeneous Calder\'on reproducing
formula and the same phenomena happens in the
Littlewood-Paley characterization
of the localized Hardy space established in \cite{hmy2}.
We also notice that when $\cx$ is an RD-space and $\mu(\cx)<\fz$,
Theorem \ref{t3.2} also holds;
moreover, in this setting, \eqref{3.1} holds automatically.

(ii) Combining Theorems \ref{t3.2} with \ref{t3.4},
we characterize $H^p_{\ell,\,\rm at}(\cx)$
with $p\in(n/(n+1),\,1]$ via
the radial maximal function and the grand maximal function.
By a reason similar to that
of Remark \ref{r3.1},
the range $p\in(n/(n+1),\,1]$ is the best possible
for which $H^p_{\ell,\,\rm at}(\cx)$ can be characterized
by the radial maximal function.

(iii) If we relax $d$ to be a quasi-metric,
then Theorems  \ref{t3.2} through \ref{t3.4}
still hold by replacing $p\in(n/(n+1),1]$
with $p\in(n/(n+\tz),\,1]$
 for any $( \tz,\,
1,\,1)-\ati$  $\{S_k\}_{k\in\zz}$,
where $\tz\in(0,\,1)$ is the same as in Remark \ref{r2.2} (iv).

(iv) Theorems \ref{t3.2} through \ref{t3.4}
are also true if we replace the discrete approximation of the identity
by the continuous one
as in Remark \ref{r2.2} (ii).
\end{rem}

\section{Applications to some differential operators}\label{s4}

Beyond the Ahlfors $n$-regular metric measure space
and Lie groups of polynomial growth,
we list several other specific settings where Theorems \ref{t3.1}
through \ref{t3.4} work.

{\bf (I)  Hardy spaces associated to a certain Laplace-Beltrami operator}

Let $M$ be a complete noncompact connected Riemannian manifold,
$d$ the geodesic distance,
$\mu$ the Riemannian measure and $\nabla$ the Riemannian gradient.
Denote by $|\cdot|$ the length in the tangent space.
One defines $\Delta$, the Laplace-Beltrami operator,
as the self-adjoint positive operator
on $L^2(M)$ by the formal integration by parts
$\langle \Delta f,\,f\rangle=\||\nabla f|\|_{L^2(M)}^2$
for all $f\in\ccc_0^\fz(M)$.
Denote by $T_t(x,\,y)$ with $t>0$ and $x,\,y\in M$
the heat kernel of $M$,
namely, the kernel of the heat semigroup $\{e^{-t\Delta}\}_{t>0}$.
One says $T_t$ satisfies the Li-Yau type estimates if
there exist  some positive constants $C_5, \wz C_5$ and $C$ such that
\begin{equation}\label{4.1}
C^{-1}\frac1{V_{\sqrt t}(x)}\exp\lf\{-\frac{d(x,y)^2}{\wz C_5t}\r\}\le
T_t(x,\,y)\le C
\frac1{V_{\sqrt t}(x)}\exp\lf\{-\frac{d(x,y)^2}{C_5t}\r\}
\end{equation}
for all $x,\,y\in M$ and $t>0$.
It is well-known \cite{ly86} that such estimates
hold on manifolds with non-negative Ricci curvature.
Later it has been proved in \cite{s95}
that the Li-Yau type estimates are equivalent to
the conjunction of the doubling property \eqref{1.1} and the scaled
Poincar\'e inequality that
for every ball $B\equiv B(x,\,r)$ and every $f$ with $f,\,\nabla f\in L^2_\loc(M)$,
 \begin{equation}\label{4.2}
 \int_B|f(y)-f_B|^2\,d\mu(y)\le Cr^2\int_B|\nabla f(y)|^2\,d\mu(y),
 \end{equation}
where $f_B$ denotes the average of $f$ on $B$, namely,
$f_B=\frac{1}{\mu(B)}\int_B f(x)\,d\mu(x)$.

Assume that $M$ satisfies the doubling property \eqref{1.1}
 and supports a scaled Poincar\'e inequality \eqref{4.2}.
 Then $(M,\,d,\,\mu)$ is a connected space of homogeneous type and hence,
 an RD-space. Observe that the Li-Yau type estimates
do  already imply some regularity estimates for the heat kernel:
there exist positive constants $C$ and $\az\in(0,\,1)$ such that
\begin{equation}\label{4.3}
|T_t(x,\,y)-T_t(z,\,y)|\le \lf(\frac{d(x,\,z)}{\sqrt t}\r)^\az
\frac{C}{V_{\sqrt t}(y)}
\end{equation}
for all $x,\,y,\,z\in M$ and $t>0$;
see, for example, \cite{s95,gr}.
See also \cite{acdh} for more discussions on the regularity of the heat kernels
and connections with the boundedness of Riesz transforms in this setting.

Notice that $e^{-t\Delta}1=1$; see \cite{g86}.
From this, \eqref{4.1}, \eqref{4.3}, the semigroup property and Remark \ref{r2.2},
it follows that
$\{T_{t^2}\}_{t>0}$ is just a continuous $(\az',\,N,\,N)$-$\ati$
for each $N\in\nn$ and $\az'\in(0,\,\az)$ as in
Definition \ref{d2.1}.
Define the semigroup maximal function and the localized one by
$T^+(f)(x)\equiv\sup_{t>0}|e^{-t\Delta}(f)(x)|$  and
$$T^+_{2^{-\ell}}(f)(x)\equiv\sup_{0<t<2^{-\ell}}
|e^{-t\Delta}(f)(x)|+\frac1{V_{2^{-\ell}}(x)}
\int_{B(x,\,2^{-\ell})}|T_{2^{-\ell}}(f)(y)|\,d\mu(y)$$
for all $\ell\in\zz$,
suitable distributions $f$ and $x\in\cx$.

Let $n$ be the same as in \eqref{1.1}. For $p\in(n/(n+\az),\,1]$,
the Hardy space $H^p_\Delta(M)$ and the localized one
$H^p_{\ell,\,\Delta}(M)$ with $\ell\in\zz$
on $M$ associated to the Laplace-Beltrami
operator $\Delta$ are defined, respectively, by
\begin{equation}\label{4.4}
H^p_\Delta(M)\equiv\lf\{f\in (\cg^\ez_0(\bz,\,\gz))' : \
 \|f\|_{ H_{\Delta}^p(M)}=\|T^+(f)\|_{L^p(\cx)}<\fz\r\},
 \end{equation}
 and
\begin{equation}\label{4.5}
H^p_{\ell,\,\Delta}(M)\equiv\lf\{f\in (\cg^\ez_0(\bz,\,\gz))' : \
 \|f\|_{ H_{\ell,\,\Delta}^p(M)}=\|T^+_{2^{-\ell}}(f)\|_{L^p(\cx)}<\fz\r\},
 \end{equation}
 where $\ez\in(n(1/p-1),\,\az)$ and $\bz,\,\gz\in(n(1/p-1),\,\ez)$.
Applying Theorems \ref{t3.1}, \ref{t3.2} and \ref{t3.4},
we have the following result.

\begin{prop}\label{p4.1} Let $p\in(n/(n+\az),\,1]$ and $\ell\in\zz$.

(i) Then $H^p_{\Delta}(M)= H_{\rm at}^p(M)$ with equivalent norms.

(ii) Then $H^p_{\ell,\,\Delta}(M)= H_{\ell,\,\rm at}^p(M)$ with
equivalent norms uniformly in $\ell$.
\end{prop}

\begin{rem}\label{r4.1}\rm (i) Proposition \ref{p4.1} also implies that
the definition of the (localized) Hardy spaces are independent of the choices
of $\ez\in(n(1/p-1),\,\az)$ and $\bz,\,\gz\in(n(1/p-1),\,\ez)$.

(ii)  Proposition \ref{p4.1} improves the result stated by
Russ \cite{r} and  Theorem 8.2 of Auscher, McIntosh and Russ \cite{amr}.

 (iii) According to the approach here, the range $p\in(n/(n+\az),\,1]$
 is the best possible to define $H^p_{\Delta}(\cx)$ and
 $H^p_{\ell,\,\Delta}(M)$,  and to establish
 their equivalences with the corresponding atomic Hardy spaces, respectively.
\end{rem}

{\bf (II) Hardy spaces for the degenerate Laplace on $\rrm$}

Let $m\ge 3$ and $\rrm$ be the $m$-dimensional Euclidean space
endowed with the Euclidean norm $|\cdot|$ and the Lebesgue measure $dx$.
Recall that a nonnegative locally integrable function $w$ is called
an $A_2(\rrm)$ weight in the sense of Muckenhoupt if
$$\sup_{B}\lf\{\frac1{|B|}\int_B w(x)\,dx\r\}
\lf\{\frac1{|B|}\int_B [w(x)]^{-1}\,dx\r\}<\fz,$$
where the supremum is taken over all the balls in $\rrm$.
Observe that if we set $w(E)\equiv\int_Ew(x)dx$
for any measurable set $E$, then
there exist positive constants $C,\,n$ and $\kz$
such that for all $x\in\rrm$, $\lz>1$ and $r>0$,
$$C^{-1}\lz^\kz w(B(x,\,r))\le w(B(x,\,\lz r))\le C\lz^nw(B(x,\,r)),$$
namely, the measure $w(x)\,dx$ satisfies \eqref{1.1}
and \eqref{1.2}.
Thus $(\rrm,\,|\cdot|,\,w(x)\,dx)$ is an RD-space.

Let $w\in A_2(\rrm)$ and
$\{a_{i,\,j}\}_{1\le i,\,j\le m}$ be a real
symmetric matrix function satisfying that
for all $x,\,\xi\in\rrm$,
$$C^{-1}w(x)|\xi|^{2}\le\sum_{1\le i,\,j\le m}a_{i,\,j}(x)\xi_i \xi_j
\le Cw(x)|\xi|^2.$$ Then the (possibly) degenerate elliptic operator
$\cl_0$ is defined by
$$\cl_0 f(x)\equiv-\frac1{w(x)}\sum_{1\le i,\,j\le m}
\partial_i(a_{i,\,j}(\cdot)\partial_j f)(x),$$
where $x\in\rrm$.
Denote by $\{T_t\}_{t>0}\equiv \{e^{-t\cl_0}\}_{t>0}$
the semigroup generated by $\cl_0$. We also denote the kernel of
$T_t$ by $T_t(x,\,y)$ for all $x,\,y\in\rrm$ and $t>0$.
Then it is known that there exist positive constants $C,\,C_6,\,\wz C_6$
and $\az\in(0,\,1]$ such that
for all $t>0$ and $x,\,y\in\rrm$,
$$
C^{-1}\frac1{V_{\sqrt t}(x)}\exp\lf\{-\frac{|x-y|^2}{\wz C_6t}\r\}\le
T_t(x,\,y)\le C
\frac1{V_{\sqrt t}(x)}\exp\lf\{-\frac{|x-y|^2}{C_6t}\r\};
$$
that for all $t>0$ and $x,\,y,\,y'\in\rrm$ with $|y-y'|<|x-y|/4$,
$$
 |T_t(x,\,y)-T_t(x,\,y')|\le C
\frac1{V_{\sqrt t}(x)}\lf(\frac{|y-y'|}{\sqrt t}\r)^\az\exp
\lf\{-\frac{|x-y|^2}{C_6t}\r\};
$$
and, moreover, that for all $t>0$ and $x,\,y\in\rn$,
$\int_\rrm T_t(x,\,y)\,w(y)\,dy=1$;
see, for example, Theorems 2.1, 2.3, 2.4 and 2.7, and Corollary 3.4 of
\cite{hs01}.
This together with the symmetry of the heat kernel and
Remark \ref{r2.2} implies that $\{T_t\}_{t>0}$
is a continuous $(\az,\,N,\,N)$-$\ati$ for
any $N\in\nn$. Thus, Theorems \ref{t3.1} through \ref{t3.4} in Section 3 also
work here for $p\in(n/(n+\az),\,1]$.
 Define the  Hardy spaces $H^p_{\cl_0}(\rrm,w)$ and the localized one
 $H^p_{\ell,\,\cl_0}(\rrm,w)$ with $\ell\in\zz$
 associated to the (possibly) degenerate Laplace as
 in \eqref{4.4} and \eqref{4.5} of (I).
 We then have the following conclusions.

\begin{prop}\label{p4.2} Let $p\in(n/(n+\az),\,1]$ and $\ell\in\zz$.

(i) Then $H^p_{\cl_0}(\rrm,w)= H_{\rm at}^p(\rrm,w)$ with equivalent norms.

(ii) Then $H^p_{\ell,\,\cl_0}(\rrm,w)= H_{\ell,\,\rm at}^p(\rrm,w)$ with
equivalent norms uniformly in $\ell$.
\end{prop}

{\bf (III) Hardy spaces associated to a certain sub-Laplace operator}

The following example deals with the differential operators on
the noncompact $C^\fz$-manifold $M$
arising as the boundary of an unbounded model
polynomial domain in $\cc^2$ introduced by Nagel and
Stein \cite{ns06} (see also \cite{ns04}).

Let $\boz\equiv\{(z,\,w)\in\cc^2:\,{\rm Im}[w]>P(z)\}$, where $P$ is a real,
subharmonic polynomial of degree $m$.
Then $M\equiv\partial \boz$ can be identified with
$\cc\times\rr\equiv\{(z,\,t):\ z\in\cc,\, t\in\rr\}$.
The basic (0,\,1) Levi vector field is then
$\overline Z=\frac{\partial}{\partial \bar z}-
i\frac{\partial P}{\partial \bar z}\frac{\partial}{\partial t}$,
and rewrite as $\overline Z=X_1+iX_2$.
The real vector fields $\{X_1,\,X_2\}$ and their commutators of
orders $\le m$ span the tangent space at each point of $M$.
Thus $M$ is of finite type $m$.

We denote by $d$ the control distance in $M$ generated by $X_1$ and $X_2$;
see \cite{nsw} for the definition of the control distance.
Let $\mu$ be the Lebesgue measure on $M$. Then,
Nagel and Stein showed that
there exist positive constants $Q\ge 4$ and $C$
such that for all $s\ge1$, $x\in\cx$ and $r>0$,
$$C^{-1}s^4V(x,\,r)\le
V(x,\,sr)\le  Cs^ QV(x,\,r),
$$
which implies that $(M,\,d,\,\mu)$ is an RD-space; see \cite{ns04,ns06}.

Moreover, denote by $X_j^\ast$ the formal adjoint of $X_j$, namely,
$\langle X_j^\ast\vz,\,\psi\rangle=\langle\vz,\,X_j\psi\rangle$
for $\vz,\,\psi\in C_c^\fz(M)$, where
$\langle\vz,\,\psi\rangle=\int_M\vz(x)\psi(x)\,dx.$
In general, $X_j^\ast=-X_j+a_j$, where $a_j\in C^\fz(M)$.
The sub-Laplacian $\cl$ on $M$ is formally given by
$ \cl\equiv X_1^\ast X_1+  X_2^\ast X_2.$
Set the heat operator $T_t\equiv e^{-t\cl}$ for $t>0$.

Nagel and Stein further established the non-Gaussian upper bound estimate
and regularity of $\{T_t\}_{t>0}$; see \cite{ns04,ns01}.
These properties further imply that
  $\{T_t\}_{t>0}$ forms a continuous
$(1,\,1,\,1)$-$\ati$ as in Definition \ref{d2.1}.
Thus Theorems \ref{t3.1} through \ref{t3.4} in Section 3 also
 work here. In particular,  we define the
 Hardy space $H^p_{\cl }(M)$ and the localized one
  $H^p_{\ell,\,\cl }(M)$ with $\ell\in\zz$
 associated to the sub-Laplace as in \eqref{4.4} and \eqref{4.5} of (I).
We then have the following conclusions.

 \begin{prop}\label{p4.3} Let $p\in(n/(n+1),\,1]$ and $\ell\in\zz$.

(i) Then $H^p_{\cl }(M)= H_{\rm at}^p(M)$ with equivalent norms.

(ii) Then $H^p_{\ell,\,\cl }(M)= H_{\ell,\,\rm at}^p(M)$ with
equivalent norms uniformly in $\ell$.
\end{prop}

\section{Proofs of main theorems}\label{s5}

We need the Calder\'on reproducing formula to prove Theorems
\ref{t3.1} through \ref{t3.4}. First we recall the dyadic cubes on
spaces of homogeneous type constructed by Christ \cite{ch}.

\begin{lem}\label{l5.1}
Let $\cx$ be a space of homogeneous type. Then there exists a
collection $\{Q_\az^k\subset\cx:\ k\in\zz,\,\az\in I_k\}$ of open
subsets, where $I_k$ is certain index set, and positive constants
$\dz\in(0,\,1)$, $C_7$ and $C_8$ such that

(i) $\mu(\cx\setminus\cup_\az Q_\az^k)=0$ for each fixed $k$ and
$Q_\az^k\cap Q_\bz^k=\emptyset$ if $\az\ne\bz$;

(ii) for any $\az,\,\bz,\,k,\,\ell$ with $\ell\ge k$, either
$Q_\bz^\ell\subset Q_\az^k$ or $Q_\bz^\ell\cap Q_\az^k=\emptyset$;

(iii) for each $(k,\,\az)$ and $\ell<k$, there exists a unique $\bz$
such that $Q_\az^k\subset Q_\bz^\ell$;

(iv) $\diam(Q_\az^k)\le C_7\dz^k$;

(v) each $Q_\az^k$ contains certain ball $B(z_\az^k,\,C_8\dz^k)$, where
$z_\az^k\in\cx$.
\end{lem}

In fact, we can think of $Q_\az^k$ as being a dyadic cube with
diameter roughly $\dz^k$  and centered at $z_\az^k$. In what follows, for
simplicity, we may assume that $\dz=1/2$; see \cite[p.\,25]{hmy2}
as to how to remove this restriction.

Fix $j_0\in \nn$ such that
$2^{-j_0}C_7<1/3.$
For any given $j\ge j_0$, and for any $k\in\zz$ and $\tau\in I_k$,
we denote by
$Q_\tau^{k,\,\nu}$, \, $\nu=1,\,2,\,\cdots,\,N(k,\,\tau)$, the set of
all cubes $Q_{\tau'}^{k+j}\subset Q_\tau^k$.
Denote by $z_\tau^{k,\,\nu}$ the center of $Q_\tau^{k,\,\nu}$ and
let
$y_\tau^{k,\,\nu}$ be a point in $Q_\tau^{k,\,\nu}$.

For any given $\ell\in\zz\cup\{-\fz\}$ and
$j\ge j_0$,  we pick a point
$y_\tau^{k,\,\nu}$ in the cube $Q_\tau^{k,\,\nu}$
for each integer $k\ge \ell$ when $\ell\in\zz$ or $k\in\zz$
when $\ell=-\fz$,\,$\tau\in I_k$ and
$\nu=1,\,\cdots,\, N(k,\,\tau)$, and denote by $\cd(\ell,\,j)$ the
set of all these points; namely,
 $\cd(-\fz,\,j)\equiv \{y_\tau^{k,\,\nu}:\
k\in\zz,\ \tau\in I_k, \ \nu=1,\,\cdots,\,
N(k,\,\tau) \}$
 and
 $\cd(\ell,\,j)\equiv\{y_\tau^{k,\,\nu}:\
k=\ell,\,\cdots,\ \tau\in I_k, \ \nu=1,\,\cdots,\,
N(k,\,\tau)\}.$

The following Calder\'on reproducing formula was established in
Theorem 4.1 of \cite{hmy2}.

\begin{thm}\label{t5.1}
Let $\ez_1\in(0,\,1]$,\,$\ez_2,\,\ez_3>0$,
$\ez\in(0,\,\ez_1\wedge\ez_2)$ and $\{S_k\}_{k\in\zz}$ be an
$(\ez_1,\,\ez_2,\,\ez_3)$-$\ati$. Set $D_k\equiv S_k-S_{k-1}$ for $k\in\zz$.
Then there exists $j>j_0$ such that for any choice of
$\cd(-\fz,\,j)$, there exists a sequence of
operators $\{\wz D_k \}_{k\in\zz}$ such that for all
$f\in(\ocg^\ez_0(\bz,\,\gz))'$ with $\bz,\,\gz\in(0,\,\ez)$, and $x\in\cx$,
\begin{eqnarray*}
f(x)&&=\sum_{k=-\fz}^\fz\sum_{\tau\in I_k}\sum_{\nu=1}^{N(k,\,\tau)}
\mu(Q_\tau^{k,\,\nu})\wz D_{k}(x,\,y^{k,\,\nu}_\tau)
D_k(f)(y^{k,\,\nu}_\tau),
\end{eqnarray*}
where the series converge  in $(\ocg^\ez_0(\bz,\,\gz))'$. Moreover,
for any $\ez'\in[\ez,\,\ez_1\wedge\ez_2)$, there exists a positive constant
$C$ depending on $\ez'$, but not on $j$ and $\cd(-\fz,\,j)$, such that
$\{\wz D_k(x,\,y)\}_{k\in\zz}$, the kernels of
$\{\wz D_k\}_{k\in\zz}$, satisfy (i) and (ii) of
Definition \ref{d2.1}
with $\ez_1$ and $\ez_2$ replaced by $\ez'$ and the constant
$C_4$ replaced by $C$, and $\int_\cx \wz
D_k(x,\,y)\,d\mu(x)=0=\int_\cx \wz D_k(x,\,y)\,d\mu(y)$.
\end{thm}

The following inhomogeneous Calder\'on reproducing formula
was established in \cite{hmy2} and \cite{gly2}. For a measurable set $E$, we
set $m_E(f)\equiv\frac1{\mu(E)}\int_E f(y)\,d\mu(y)$.

\begin{thm}\label{t5.2}
Let $\ez_1\in(0,\,1]$,\,$\ez_2,\,\ez_3>0$,
$\ez\in(0,\,\ez_1\wedge\ez_2)$ and $\{S_k\}_{k\in\zz}$ be an
$(\ez_1,\,\ez_2,\,\ez_3)$-$\ati$.  Set $D_k\equiv S_k-S_{k-1}$ for $k\in\zz$.
Then there exist  $N\in\nn$ and $j_1>j_0$
such that for all $j>j_1$,
 $\ell\in\zz$ and $\cd(\ell+N+1,\,j)$,
there exist operators $\{\wz D_k(x,\,y)\}_{k=\ell}^\fz$ such
that for any $f\in(\cg^\ez_0(\bz,\,\gz))'$ with
$\bz,\,\gz\in(0,\,\ez)$, and all $x\in\cx$,
\begin{eqnarray*}
f(x)&&= \sum_{\tau\in I_\ell}\sum_{\nu=1}^{N(\ell,\,\tau)}
\int_{Q_\tau^{\ell,\,\nu}}\wz D_\ell(x,\,y)d\mu(y)\,
m_{Q^{\ell,\,\nu}_\tau}(S_\ell(f))\\
&&\quad+ \sum_{k=\ell+1}^{\ell+N} \sum_{\tau\in
I_k}\sum_{\nu=1}^{N(k,\,\tau)}
\int_{Q_\tau^{k,\,\nu}}\wz D_k(x,\,y)d\mu(y)\, m_{Q^{k,\,\nu}_\tau}(D_k(f))\\
&&\quad+\sum_{k=\ell+N+1}^\fz\sum_{\tau\in
I_k}\sum_{\nu=1}^{N(k,\,\tau)} \mu(Q_\tau^{k,\,\nu})\wz
D_{k}(x,\,y^{k,\,\nu}_\tau) D_k(f)(y^{k,\,\nu}_\tau),
\end{eqnarray*}
where the series converge  in $(\cg^\ez_0(\bz,\,\gz))'$. Moreover,
for any $\ez'\in[\ez,\,\ez_1\wedge\ez_2)$, there exists a positive constant $
C$ depending on $N$ and $\ez'$, but not on $\ell $, $j$ and
$\cd( \ell+N+1,\,j)$, such that $\wz D_k$ for $k\ge\ell$ satisfies
(i) and (ii) of Definition \ref{d2.1} with $\ez_1$ and $\ez_2$
replaced by $\ez'$ and the constant $C_4$ replaced by $ C $, and
$\int_\cx \wz D_k(x,\,y)\,d\mu(x)=\int_\cx \wz
D_k(x,\,y)\,d\mu(y)=1$ when $\ell\le k\le\ell+N$; $=0$ when
$k\ge\ell+N+1$.
\end{thm}

\begin{rem}\label{r5.1}\rm
(i) We remark that by checking the proofs of Theorem 4.1 and Theorem 4.3
in \cite{hmy2}, it is easy to see that the constants $C$
in  Theorems \ref{t5.1} and  \ref{t5.2} are independent of $j$.
This observation plays a key role in the proof of Theorem \ref{t3.2}
(ii).

(ii) For simplicity, in Theorem \ref{t5.2}, we always
assume that $N=0$ in what follows.
\end{rem}

The following technical lemma is also used;
see \cite{hmy2} for its proof.

\begin{lem}\label{l5.2}
Let $\ez\in(0,\,1]$, $ r\in(n/(n+\ez),1]$ and $j\ge j_0$.
Then there exists a positive
constant $C$ independent of $j$ such that for all $k,\,k'\in\zz$,
$a^{k,\,\nu}_\tau\in\cc$, $y^{k,\,\nu}_\tau\in Q^{k,\,\nu}_\tau$
with $\tau\in I_k$ and $\nu=1,\,\cdots,\,N(k,\,\tau)$ and $x\in\cx$,
\begin{eqnarray*}
&&\sum_{\tau\in I_k}\sum_{\nu=1}^{N(k,\,\tau)}\mu(Q^{k,\,\nu}_\tau)
 \frac{|a^{k,\,\nu}_\tau|}{V_{2^{-(k'\wedge k)}}(x)+V(x,\,y^{k,\,\nu}_\tau)}
\lf[\frac{2^{-(k'\wedge k)}}{2^{-(k'\wedge k)}+d(x,\,y^{k,\,\nu}_\tau)}\r]^{\ez}\\
&&\quad\le C2^{jn(1/r-1)}2^{[(k'\wedge k)-k]n(1-1/r)}
\lf\{\hl\lf(\sum_{\tau\in I_k}\sum_{\nu=1}^{N(k,\,\tau)}
|a^{k,\,\nu}_\tau|^r\chi_{Q^{k,\,\nu}_\tau}\r)(x)\r\}^{1/r}.
\end{eqnarray*}
\end{lem}

Now we prove Theorem \ref{3.1} by using the homogeneous
Calder\'on formula given in Theorem \ref{t5.1}.

\newtheorem{pfone}{\it Proof of Theorem \ref{t3.1}}
\renewcommand\thepfone{}

\begin{pfone}\rm To prove Theorem \ref{t3.1}, it suffices to prove that for any
fixed $\bz,\,\gz\in(n(1/p-1),\,\ez)$ ,
there exist positive constants $r\in (n/(n+\bz\wedge\gz),\,p)$ and $C$
such that for all $f\in(\cg_0^\ez(\bz,\,\gz))'$
and $x\in\cx$,
\begin{equation}\label{5.1}
G(f)(x)\le C\lf\{\hl([S^+(f)]^r)(x)\r\}^{1/r}.
\end{equation}
To this end, it suffices to prove that for any $\vz\in\ocg^\ez_0(\bz,\,\gz)$
with $\|\vz\|_{\cg^\ez_0(x,\,2^{-\ell},\,\bz,\,\gz)}\le1$,
\begin{equation}\label{5.2}
|\langle f,\,\vz\rangle| \le C\lf\{\hl([S^+(f)]^r)(x)\r\}^{1/r}.
\end{equation}

Assume that \eqref{5.2} holds for the moment.
Then for any $\phi\in\cg^\ez_0(\bz,\,\gz)$
with $$\|\phi\|_{\cg^\ez_0(x,\,2^{-\ell},\,\bz,\,\gz)}\le1$$
and  $\sz\equiv\int_\cx \phi(y)\,d\mu(y)$,
we have
$$|\sz|
\le \int_\cx \frac1{V_{2^{-\ell}}(x)+V(x,\,y)}
\lf(\frac{2^{-\ell}}{2^{-\ell}+d(x,\,y)}\r)^\gz\,d\mu(y)\ls1.$$
Set $\vz(y)\equiv\frac1{1+|\sz| C_3}[\phi(y)-\sz S_\ell(x,\,y)]$.
Then $\int_\cx \vz(y)\,d\mu(y)=0$
and hence, $\vz\in\ocg^\ez_0(\bz,\,\gz)$
with $\|\vz\|_{\cg^\ez_0(x,\,2^{-\ell},\,\bz,\,\gz)}\le 1$.
Since
$$|\langle f,\,\phi\rangle|
\le |\sz||S_\ell(f)(x)|+(1+|\sz|C_3)|\langle f,\,\vz\rangle|,$$ by
\eqref{5.2} and taking the supremum over all $\ell\in\zz$ and
$\phi\in\cg^\ez_0(\bz,\,\gz)$ with
$\|\phi\|_{\cg^\ez_0(x,\,2^{-\ell},\,\bz,\,\gz)}\le1$, we further obtain that
for all $x\in\cx $,
$$G(f)(x)\le |\sz| S^+(f)(x)+ C(1+\sz C_3)\lf\{\hl([S^+(f)]^r)(x)\r\}^{1/r},$$
which yields \eqref{5.1}.

To prove \eqref{5.2}, for each $Q^{k,\,\nu}_\tau$, we pick a point
$y^{k,\,\nu}_\tau\in Q^{k,\,\nu}_\tau$ such that
$$|D_k(f)(y^{k,\,\nu}_\tau)|\le2\inf_{z\in Q^{k,\,\nu}_\tau}|D_k(f)(z)|,$$
which further implies that
$$|D_k(f)(y^{k,\,\nu}_\tau)|
\le  2\inf_{z\in Q^{k,\,\nu}_\tau}(|S_k(f)(z)|+ |S_{k-1}(f)(z)|)
\le 4\inf_{z\in Q^{k,\,\nu}_\tau}S^+(f)(z).
$$
Denote by $\cd(-\fz,\,j)$ the collection of all such points $y^{k,\,\nu}_\tau$.
Observe that $\ocg^\ez_0(\bz,\,\gz)$ is a subspace of
$\cg^\ez_0(\bz,\,\gz)$ and for any
$\vz\in\ocg^\ez_0(\bz,\,\gz)$,
$\|\vz\|_{\ocg^\ez_0(\bz,\,\gz)}
=\|\vz\|_{\cg^\ez_0(\bz,\,\gz)}$.
So for any
$f\in(\cg^\ez_0(\bz,\,\gz))'$,
we know that $f$ uniquely induces a bounded linear functional $g$ on
$\ocg^\ez_0(\bz,\,\gz)$, namely, $g\in(\ocg^\ez_0(\bz,\,\gz))'$;
moreover, $\|g\|_{(\ocg^\ez_0(\bz,\,\gz))'}
\le \|f\|_{(\cg^\ez_0(\bz,\,\gz))'}$.
By $\vz\in \ocg^\ez_0(\bz,\,\gz)$, $D_k(y,\,\cdot)\in \ocg^\ez_0(\bz,\,\gz)$
for all $y\in\cx$ and Theorem \ref{t5.1}, we have
\begin{eqnarray*}
\langle f,\,\vz\rangle\equiv\langle g,\,\vz\rangle&&=\sum_{k=-\fz}^\fz
\sum_{\tau\in I_k}\sum_{\nu=1}^{N(k,\,\tau)}
\mu(Q_\tau^{k,\,\nu})\wz D^\ast_{k}(\vz)(y^{k,\,\nu}_\tau)
D_k(g)(y^{k,\,\nu}_\tau)\\
&&=\sum_{k=-\fz}^\fz
\sum_{\tau\in I_k}\sum_{\nu=1}^{N(k,\,\tau)}
\mu(Q_\tau^{k,\,\nu})\wz D^\ast_{k}(\vz)(y^{k,\,\nu}_\tau)
D_k(f)(y^{k,\,\nu}_\tau),
\end{eqnarray*}
where and in what follows, for all $y\in\cx$,
$$\wz D^\ast_{k}(\vz)(y)=\int_\cx \wz D_k(z,\,y)\vz(z)\,d\mu(z).$$

We recall that there exist positive constants $\bz'\in(n(1/p-1),\,\bz)$ and
$C$ such that for all $y\in\cx$,
\begin{equation}\label{5.3}
|\wz D^\ast_{k}(\vz)(y)|\le C2^{-|k-\ell|\bz'}
\frac1{V_{2^{-(\ell\wedge k)}}(x)+V(x,\,y)}\lf(\frac{2^{-(\ell\wedge k)}}
{2^{-(\ell\wedge k)}+d(x,\,y)}\r)^{\gz};
\end{equation}
see the proof of Proposition 5.7 of \cite{hmy2} for the details.
Recall that $\ell\wedge k\equiv\min\{\ell,\,k\}$. We point out that,
to obtain the decay factor $2^{-|k-\ell|\bz'}$ in \eqref{5.3},
we need to use the fact that
$\int_\cx \vz(z)\,d\mu(z)=0$ and
$\int_\cx \wz D_k(z,\,y)\,d\mu(z)=0$ for all $y\in\cx$.

Then by Lemma \ref{l5.2} with $r\in(n/(n+\bz'\wedge\gz),\,p)$, we have
\begin{eqnarray*}
|\langle f,\,\vz\rangle|&&\ls\sum_{k=-\fz}^\fz 2^{-|k-\ell|\bz'}
\sum_{\tau\in I_k}\sum_{\nu=1}^{N(k,\,\tau)}
\mu(Q_\tau^{k,\,\nu})\inf_{z\in Q^{k,\,\nu}_\tau}S^+(f)(z)\\
&&\quad \times
\frac1{V_{2^{-(\ell\wedge k)}}(x)+V(x,\,y^{k,\,\nu}_\tau)}
\lf(\frac{2^{-(\ell\wedge k)}}
{2^{-(\ell\wedge k)}+d(x,\,y^{k,\,\nu}_\tau)}\r)^{\gz}\\
&&\ls\sum_{k=-\fz}^\fz 2^{-|k-\ell|\bz'}
2^{[(\ell\wedge k)-k]n(1-1/r)}\\
 &&\quad\times
\lf\{\hl\lf(\sum_{\tau\in I_k}\sum_{\nu=1}^{N(k,\,\tau)}
\lf|\inf_{z\in Q^{k,\,\nu}_\tau}S^+(f)(z)\r|^r
\chi_{Q^{k,\,\nu}_\tau}\r)(x)\r\}^{1/r}\\
&&\ls \lf( \sum_{k=-\fz}^\ell 2^{-(\ell-k)\bz'}+
\sum_{k=\ell+1}^\fz 2^{-(k-\ell)[\bz'-n(1/r-1)]}\r)
\lf\{\hl\lf([ S^+(f)]^r\r)(x)\r\}^{1/r}\\
&&\ls \lf\{\hl\lf([ S^+(f)]^r\r)(x)\r\}^{1/r},
\end{eqnarray*}
which gives \eqref{5.2}.
This finishes the proof of Theorem \ref{t3.1}.

\end{pfone}

To prove Theorem \ref{t3.2},  we first introduce a variant of the localized
radial maximal function.
\begin{defn}\label{d5.1}
Let $\ez_1\in(0,\,1]$,\,$\ez_2,\,\ez_3>0$,
$\ez\in(0,\,\ez_1\wedge\ez_2)$ and $\{S_k\}_{k\in\zz}$ be an
$(\ez_1,\,\ez_2,\,\ez_3)$-$\ati$. Let
$p\in(n/(n+\ez),\,1]$ and $\bz,\,\gz\in(n(1/p-1),\,\ez)$.
For $j\ge j_0$ and $f\in(\cg^\ez_0(\bz,\,\gz))'$, define
$$S^{+,\,j}_\ell(f)\equiv
 S^{+}_\ell(f) + \lf\{
 \sum_{\tau\in I_{\ell}}\sum_{\nu=1}^{N({\ell},\,\tau)}
[m_{Q^{{\ell},\,\nu}_\tau}(|S_{\ell}(f)|)]^p\chi_{Q^{\ell,\,\nu}_\tau}\r\}^{1/p},$$
where $\{Q^{\ell,\,\nu}_\tau\}_{\nu=1}^{N(\ell,\,\tau)}$
is the collection of all dyadic cubes  $Q^{\ell+j}_{\tau'}$
contained in  $Q^{\ell}_{\tau }$.
\end{defn}

\newtheorem{pfsecond}{\it Proof of Theorem \ref{t3.2}}
\renewcommand\thepfsecond{}

\begin{pfsecond}\rm To prove (i), we first claim that
for $p\in(n/(n+\ez),1]$, all $\ell\in\zz$ and $f\in(\cg^0_\ez(\bz,\,\gz))'$,
\begin{equation}\label{5.4}
\|S^{(a)}_\ell(f)\|^p_{L^p(\cx)}\sim
\sum_{\tau\in I_\ell}\sum_{\nu=1}^{N(\ell,\,\tau)}
[m_{Q^{\ell,\,\nu}_\tau}(|S_{\ell}(f)|)]^p\mu(Q^{\ell,\,\nu}_\tau),
\end{equation}
where the constants depend on $j,\,p$ and $a$, but not on $\ell$ and $f$.

To see this, we observe that for each given $j\ge0$,
there exists a positive constant $C(j)$
depending on $j$ such that for all $\ell\in\zz$ and $\tau\in I_\ell$,
$N(\ell,\,\tau)\le C(j)$. In fact, by Lemma \ref{l5.1},
\eqref{1.1} and \eqref{1.2},
we have
\begin{eqnarray*}
\mu(Q^{\ell}_\tau)=\sum_{\nu=1}^{N(\ell,\,\tau)}\mu(Q^{\ell,\,\nu}_\tau)
&&\ge \sum_{\nu=1}^{N(\ell,\,\tau)}V_{C_82^{-\ell-j}}(z^{\ell,\,\nu}_\tau)
 \gs \sum_{\nu=1}^{N(\ell,\,\tau)}V_{2C_72^{-\ell}}(z^{\ell,\,\nu}_\tau)
\gs  N(\ell,\,\tau)\mu(Q^{\ell}_\tau),
\end{eqnarray*}
which implies our claim.
Then it is easy to see that for all $j\ge 0$,
\begin{equation}\label{5.5}\sum_{\tau\in I_\ell}\sum_{\nu=1}^{N(\ell,\,\tau)}
[m_{Q^{\ell,\,\nu}_\tau}(|S_{\ell}(f)|)]^p\mu(Q^{\ell,\,\nu}_\tau)\sim
\sum_{\tau\in I_\ell}
[m_{Q^{\ell}_\tau}(|S_{\ell}(f)|)]^p\mu(Q^{\ell}_\tau),
\end{equation}
where the equivalent constants are independent of $f$ and $\ell$.
For any given $a>0$, let $j$ be large enough such that
$a> C_72^{-\ell-j}\ge \diam(Q^{\ell,\,\nu}_\tau)$.
Since for any $x\in Q^{\ell,\,\nu}_\tau$,
$$m_{Q^{\ell,\,\nu}_\tau}(|S_{\ell}(f)|)
\le \frac {V_{a2^{-\ell}}(x)}{\mu(Q^{\ell,\,\nu}_\tau)}
S^{(a)}_\ell(f)(x)\ls \frac
{V_{a2^{-\ell+1}}(z^{\ell,\,\nu}_\tau)}{
V_{C_82^{-\ell}}(z^{\ell,\,\nu}_\tau)}S^{(a)}_\ell(f)(x)\ls
S^{(a)}_\ell(f)(x),$$ we know that the left hand side of
\eqref{5.4} controls its right hand side
for this $j$ and thus for all
$j\ge0$ by \eqref{5.5}. For $a\ge 0$ and each $Q^\ell_\tau$, by \eqref{1.2} and a
similar argument, we have
\begin{eqnarray*}\int_{Q^\ell_\tau}[S^{(a)}_\ell(f)(x)]^p\,d\mu(x)
&&\ls
[\mu(Q^\ell_\tau)]^{1-p} \lf\{\int_{B(x^\ell_\tau,\,(a+C_7)2^{-\ell})}
S_\ell(f)(y)\,d\mu(y)\r\}^{p}\\
&&\ls \sum_{Q^\ell_{\tau'}\cap B(x^\ell_\tau,\,(a+C_7)2^{-\ell})\ne\emptyset}
[m_{Q^\ell_{\tau'}}(|S_\ell(f)|)]^p\mu(Q^\ell_{\tau'}).
\end{eqnarray*}
Similarly to above, using Lemma \ref{5.5}, \eqref{1.1} and \eqref{1.2},  we have
that $\sharp\{\tau':\
Q^\ell_{\tau'}\cap B(x^\ell_\tau,\,(a+C_5)2^{-\ell})\ne\emptyset\}$
is  bounded uniformly in $\ell$ and $\tau$,
from which it follows that the right hand side of \eqref{5.4}
controls its left hand side for $j=0$
and thus for all $j$ by \eqref{5.5}.
Here, $\sharp E$ denotes the cardinality of the set $E$.

By Lemma \ref{l5.1}  and the H\"older inequality,
we have
$$\|S_\ell(f)\|^p_{L^p(\cx)}=\sum_{\tau\in I_\ell}\sum_{\nu=1}^{N(\ell,\,\tau)}
m_{Q^{\ell,\,\nu}_\tau}(|S_{\ell}(f)|^p)\mu(Q^{\ell,\,\nu}_\tau)
\le \sum_{\tau\in I_\ell}\sum_{\nu=1}^{N(\ell,\,\tau)}
[m_{Q^{\ell,\,\nu}_\tau}(|S_{\ell}(f)|)]^p\mu(Q^{\ell,\,\nu}_\tau),$$
which together with \eqref{5.4} leads to
$$\|S^{+,j}_\ell(f)\|_{L^p(\cx)}\sim \|S^+_{\ell+1}(f)\|_{L^p(\cx)}+
\|S^{(a)}_{\ell}(f)\|_{L^p(\cx)}$$ with constants independent of
$\ell$ and $f$. Then the proof of Theorem \ref{t3.2} (i)
is reduced to showing that
\begin{equation}\label{5.6}
C^{-1}\|S^{+,j}_\ell(f)\|_{L^p(\cx)}\le
\|G_\ell(f)\|_{L^p(\cx)}\le C2^{jn(1/p-1)}\|S^{+,j}_{\ell}(f)\|_{L^p(\cx)},
\end{equation}
where $C$ is a positive constant independent of $j$, $\ell$ and $f$.

Observe that for all $x,\,y\in Q^{k,\,\nu}_\tau$ and $k\ge \ell$,
$S_\ell(f)(x)\ls G_\ell(f)(y)$, which implies the first inequality in
\eqref{5.6}.
To prove the second inequality in \eqref{5.6}, it suffices
to prove that for any fixed
$\bz,\,\gz\in(n(1/p-1),\,\ez)$, there exist positive constants $r\in
(n/(n+\bz\wedge\gz),\,p)$ and $C$
 such that for all $j\ge j_0$, $\ell\in\zz$,
 $f\in(\cg_0^\ez(\bz,\,\gz))'$ and $x\in\cx$,
\begin{equation}\label{5.7}
G_\ell(f)(x)\le C2^{jn(1/p-1)}\lf\{\hl([S^{+,\,j}_{\ell}(f)]^r)(x)\r\}^{1/r}.
\end{equation}

To prove \eqref{5.7},
we only need to prove that for all $k'\ge \ell$ and
$\vz\in\ocg^\ez_0(\bz,\,\gz)$
with $\|\vz\|_{\cg^\ez_0(x,\,2^{-k'},\,\bz,\,\gz)}\le1$,
\begin{equation}\label{5.8}
|\langle f,\,\vz\rangle| \le C2^{jn(1/p-1)}
\lf\{\hl([S^{+,\,j}_{\ell}(f)]^r)(x)\r\}^{1/r}.
\end{equation}
In fact, assume that  $\phi\in\cg^\ez_0(\bz,\,\gz)$
with $\|\phi\|_{\cg^\ez_0(x,\,2^{-k'},\,\bz,\,\gz)}\le1$ for certain $k\ge \ell$.
Set
$\sz\equiv\int_\cx \phi(y)\,d\mu(y) $
and $\vz(y)=\frac1{1+|\sz| C_3}[\phi(y)-\sz S_{k'}(x,\,y)]$.
Then we have $\int_\cx \vz(y)\,d\mu(y)=0$
and hence, $\vz\in\ocg^\ez_0(\bz,\,\gz)$
with $\|\vz\|_{\cg^\ez_0(x,\,2^{-k'},\,\bz,\,\gz)}\le 1$.
Moreover, we have
$$|\langle f,\,\phi\rangle|
\le |\sz||S_{k'}(f)(x)|+(1+|\sz| C_3)|\langle f,\,\vz\rangle|.$$
By taking the supremum over
all $k'\ge\ell$ and $\phi\in\cg^\ez_0(\bz,\,\gz)$ with
$\|\phi\|_{\cg^\ez_0(x,\,2^{-k'},\,\bz,\,\gz)}\le1$ together with \eqref{5.8},
we have
$$G_{\ell}(f)(x)\le
\sz S^{+,\,j}_{\ell}(f)(x)+ C(1+\sz
C_3)2^{jn(1/p-1)}\lf\{\hl([S^{+,\,j}_{\ell}(f)]^r)(x)\r\}^{1/r},$$
which combined with the uniform boundedness of  $\sz$
yields \eqref{5.7}.

To prove \eqref{5.8}, for each $Q^{k,\,\nu}_\tau$ with $k\ge \ell+1$,
 we choose a point $y^{k,\,\nu}_\tau\in Q^{k,\,\nu}_\tau$
such that
$$|D_k(f)(y^{k,\,\nu}_\tau)|\le 2\inf_{z\in Q^{k,\,\nu}_\tau}
|D_k(f)(z)|\le
4\inf_{z\in Q^{k,\,\nu}_\tau}
S^{+,\,j}_{\ell}(f)(z)$$
and denote by $\cd(\ell+1,\,j)$ the collection of all such points $y^{k,\,\nu}_\tau$.
Then by Theorem \ref{t5.2}, we write
\begin{eqnarray*}
\langle f,\,\vz\rangle&&=
\sum_{\tau\in I_\ell}\sum_{\nu=1}^{N(\ell,\,\tau)}
\int_{Q_\tau^{k,\,\nu}} \wz D^\ast_{\ell}(\vz)(z)\,d\mu(z)
m_{Q_\tau^{\ell,\,\nu}}(S_\ell(f))\\
 &&\hs+
\sum_{k=\ell+1}^\fz
\sum_{\tau\in I_k}\sum_{\nu=1}^{N(k,\,\tau)}
\mu(Q_\tau^{k,\,\nu})\wz D^\ast_{k}(\vz)(y^{k,\,\nu}_\tau)
D_k(f)(y^{k,\,\nu}_\tau).
\end{eqnarray*}

Similarly to the proof of \eqref{5.3}, there exist positive constants
$\bz'\in(n(1/p-1),\,\bz)$ and $C$, independent of $j$, $\ell$ and
$\cd(\ell+1,\,j)$, such that for all $k\in\zz$ and $y\in\cx$,
$$
|\wz D^\ast_{k}(\vz)(y)|\le C2^{-|k-k'|\bz'}
\frac1{V_{2^{-(k'\wedge k)}}(x)+V(x,\,y)}\lf(\frac{2^{-(k'\wedge k)}}
{2^{-(k'\wedge k)}+d(x,\,y)}\r)^{\gz}.
$$
Then by Lemma \ref{l5.2} with $r\in(n/(n+\bz'\wedge\gz),\,p)$, we have
\begin{eqnarray*}
|\langle f,\,\vz\rangle|&&\ls
\sum_{k=\ell}^\fz 2^{-|k-k'|\bz'}
\sum_{\tau\in I_k}\sum_{\nu=1}^{N(k,\,\tau)}
\mu(Q_\tau^{k,\,\nu})\inf_{z\in Q^{k,\,\nu}_\tau}
S^{+,\,j}_{\ell}(f)(z)\\
&&\quad\times
\frac1{V_{2^{-(k'\wedge k)}}(x)+V(x,\,y^{k,\,\nu}_\tau)}
\lf(\frac{2^{-(k'\wedge k)}}
{2^{-(k'\wedge k)}+d(x,\,y^{k,\,\nu}_\tau)}\r)^{\gz}\\
&&\ls\sum_{k=\ell}^\fz 2^{-|k-k'|\bz'}
2^{[(k'\wedge k)-k]n(1-1/r)}2^{jn(1/r-1)}\\
 &&\quad\times
\lf\{\hl\lf(\sum_{\tau\in I_k}\sum_{\nu=1}^{N(k,\,\tau)}
\lf|\inf_{z\in Q^{k,\,\nu}_\tau}S^{+,\,j}_{\ell}(f)(z)\r|^r
\chi_{Q^{k,\,\nu}_\tau}\r)(x)\r\}^{1/r}\\
&&\ls 2^{jn(1/r-1)}\lf\{\hl\lf([ S^{+,\,j}_{\ell}(f)]^r\r)(x)\r\}^{1/r}.
\end{eqnarray*}
 Thus \eqref{5.8} holds, which gives (i).

To prove (ii), by the fact that $S_\ell^+(f)(y)\le G_\ell(f)(y)$
for all $y\in\cx$, we only need to prove that
for any $\bz,\,\gz\in(n(1/p-1),\,\ez)$,
there exist $r\in(n/(n+\ez),\,p)$ and a
positive constant $C$ such that for all $\ell\in\zz$,
 $j>j_0$ and $f\in(\cg^\ez_0(\bz,\,\gz))'$,
\begin{equation}\label{5.9}
G_\ell(f)\le C2^{jn(1/r-1)}\lf\{\hl([S^+_\ell(f)]^r)\r\}^{1/r}+
C2^{-j(\ez_1+n-n/r)}\lf\{\hl([G_\ell(f)]^r)\r\}^{1/r}.
\end{equation}

Assume that \eqref{5.9} holds for the moment.
Then, by the $L^{p/r}(\cx)$-boundedness of
$\hl$, we have
\begin{eqnarray*}
\|G_\ell(f)\|_{L^p(\cx)}
\le C2^{jn(1/r-1)}\|S^+_\ell(f)\|_{L^p(\cx)}
+C2^{-j(\ez_1+n-n/r)}\|G_\ell(f)\|_{L^p(\cx)}.
\end{eqnarray*}
Notice that from the assumption
\eqref{3.1} and Theorem \ref{t3.2} (i), we deduce that
$\|G_\ell(f)\|_{L^p(\cx)}<\fz$.
Thus, choosing $j$ large enough such that
$C2^{-j(\ez_1+n-n/r)}\le1/2$, we obtain $$\|G_\ell(f)\|_{L^p(\cx)}
\le C2^{jn(1/r-1)}\|S^+_\ell(f)\|_{L^p(\cx)}.$$

To prove \eqref{5.9}, for each $Q^{k,\,\nu}_\tau$ with $k\ge\ell$,
we pick a point $y^{k,\,\nu}_\tau\in Q^{k,\,\nu}_\tau$ such that
$$|S_\ell(f)(y^{\ell,\,\nu}_\tau)|
\le 2\inf_{z\in Q^{\ell,\,\nu}_\tau}S^+_{\ell}(f)(z)$$
and for $k\ge\ell+1$,
$$|D_k(f)(y^{k,\,\nu}_\tau)|
\le 2\inf_{z\in Q^{k,\,\nu}_\tau}|D_k(f)(z)
|\le 4\inf_{z\in Q^{k,\,\nu}_\tau}S^+_{\ell}(f)(z).$$
Denote by $\cd(\ell,\,j)$ the collection of all such points $y^{k,\,\nu}_\tau$.
For any $k'\ge \ell$  and any $\vz\in\ocg^\ez_0(\bz,\,\gz)$
with $\|\vz\|_{\ocg^\ez_0(x,\,2^{-k'},\,\bz,\,\gz)}\le1$, by Theorem \ref{t5.2},
 we have
\begin{eqnarray*}
\langle f,\,\vz\rangle&&
=\sum_{\tau\in I_\ell}
\sum_{\nu=1}^{N(\ell,\,\tau)}\int_{Q_\tau^{\ell,\,\nu}}
D^\ast_\ell(\vz)(y)\,d\mu(y)
[m_{Q^{\ell,\,\nu}_\tau}(S_\ell(f))-S_\ell(f)(y^{\ell,\,\nu}_\tau)]\\
&&\quad+\sum_{\tau\in I_\ell}
\sum_{\nu=1}^{N(\ell,\,\tau)}\int_{Q_\tau^{\ell,\,\nu}}
D^\ast_\ell(\vz)(y)\,d\mu(y)
 S_\ell(f)(y^{\ell,\,\nu}_\tau) \\
&&\quad+\sum_{k=\ell+1}^\fz
\sum_{\tau\in I_k}\sum_{\nu=1}^{N(k,\,\tau)}
\mu(Q_\tau^{k,\,\nu})\wz D^\ast_{k}(\vz)(y^{k,\,\nu}_\tau)
D_k(f)(y^{k,\,\nu}_\tau)\equiv J_1+J_2+J_3.
\end{eqnarray*}

Similarly to the proofs of \eqref{5.8} and \eqref{5.2},  we have
\begin{eqnarray*}
&&|J_2|+|J_3|
\ls2^{jn(1/r-1)}\lf\{\hl\lf([ S^+_\ell(f)]^r\r)(x)\r\}^{1/r}.
\end{eqnarray*}

To estimate $J_1$, observe that
there exists a positive constant $C$, independent of
$j$ and $\ell$, such that
for all
$y,\,z,\,y^{\ell,\,\nu}_\tau\in Q^{\ell,\,\nu}_\tau$,
$$\|S _\ell(y,\,\cdot)-S_\ell(y^{\ell,\,\nu}_\tau,\,\cdot)\|
_{\cg^\ez_0(z,\,
2^{-\ell},\,\bz,\,\gz)}\ls2^{-j\ez_1},$$
which implies that
\begin{eqnarray*}
|m_{Q^{\ell,\,\nu}_{\tau}}(S_\ell(f))-S_\ell(f)(y^{\ell,\,\nu}_\tau)|
\ls 2^{-j\ez_1}\inf_{z\in Q^{\ell,\,\nu}_\tau}
G_\ell(f)(z).
\end{eqnarray*}

Therefore, by Lemma \ref{l5.2}
and the fact that
$2^{-\ell}+d(x,\,y)\sim 2^{-\ell}+d(x,\,y^{\ell,\,\tau}_\nu)$
for all $y\in Q^{\ell,\,\nu}_\tau$,
we have
\begin{eqnarray*}
J_1
&&\ls2^{-j\ez_1}
\sum_{\tau\in I_\ell}\sum_{\nu=1}^{N(\ell,\,\tau)}
\mu(Q^{\ell,\,\nu}_\tau)\inf_{z\in Q^{\ell,\,\nu}_\tau} G_\ell(f)(z)\\
&&\quad\quad\times\frac1{V_{2^{-\ell}}(x)
+ V(x,\,y^{\ell,\,\tau}_\nu)}
\lf(\frac{2^{-\ell}}{2^{-\ell}+d(x,\,y^{\ell,\,\tau}_\nu)}\r)^\gz
\\
&&\ls 2^{-j(\ez_1+n-n/r)}\\
&&\quad\quad\times\lf\{\hl\lf(\sum_{\tau\in I_\ell}\sum_{\nu=1}^{N(\ell,\,\tau)}
\lf[\inf_{z\in Q^{\ell,\,\nu}_\tau} G _\ell(f)(z)\r]^r
\chi_{Q^{\ell,\,\nu}_\tau}\r)(x)\r\}^{1/r}\\
&&\ls2^{-j(\ez_1+n-n/r)}\lf\{\hl\lf(
\lf[G_{ \ell}(f)\r]^r\r)(x)\r\}^{1/r}.
\end{eqnarray*}
This gives \eqref{5.9} and hence, finishes the proof of Theorem
\ref{t3.2}.
\end{pfsecond}

Now we give the proof of Theorem \ref{t3.3}.

\newtheorem{pffive}{\it Proof of Theorem \ref{t3.3}}
\renewcommand\thepffive{}

\begin{pffive}\rm
Assume that $\|G_\ell(f)\|_{L^p(\cx)}<\fz$ and
$\{S_k\}_{k\in\zz}$ is an $(\ez_1,\,\ez_2,\,\ez_3)$-$\ati$.
By Theorem \ref{t3.1},
it suffices to prove that
$$\|S^+(f-S_\ell(f))\|_{L^p(\cx)}\ls\|G_\ell(f)\|_{L^p(\cx)} $$
with a constant independent of $f$ and $\ell$.
Write
\begin{eqnarray*}
\|S^+(f-S_\ell(f))\|_{L^p(\cx)}
&&\sim\|S^+_\ell(f)\|_{L^p(\cx)}
+\lf\|\sup_{k\ge\ell+1} |S_kS_{\ell}(f)|\r\|_{L^p(\cx)}\\
&&\quad\quad+\lf\|\sup_{k\le\ell} |S_k(f-S_{\ell}(f))|\r\|_{L^p(\cx)}
\equiv I_1+I_2+I_3.
\end{eqnarray*}
Obviously, $I_1\ls\|G_\ell(f)\|_{L^p(\cx)}$.
Write $S_kS_{\ell}(y,\,z)$ for the kernel of
$S_kS_{\ell}$. Let $\ez\in(0,\,\ez_1\wedge \ez_2)$ and $\bz,\gz\in(n(1/p-1),\,\ez)$.
Then
for any $x\in\cx$ and $k\ge\ell+1$, by \cite[Lemma 3.8]{hmy2},
 we have that
 $S_kS_\ell(x,\,\cdot)\in\cg_0^\ez(\bz,\,\gz)$ and
\begin{equation*}
\|S_kS_\ell(x,\,\cdot)\|_{\cg(x,\,\,2^{-\ell-1},\,\ez,\,\ez)}\ls1
\end{equation*}
with a constant independent of $\ell,\,k$ and $f$. From this, it follows
that $I_2\ls\|G_\ell(f)\|_{L^p(\cx)}.$

To estimate $I_3$, using Theorem \ref{t5.2}, for $k\le \ell$, we
have
\begin{eqnarray*}
&&S_k(I-S_\ell)f(x)\\
&&\quad= \sum_{\tau\in I_\ell}\sum_{\nu=1}^{N(\ell,\,\tau)}
\int_{Q_\tau^{\ell,\,\nu}}[S_k(I-S_\ell)\wz D_\ell](x,\,u)d\mu(u)\,
m_{Q^{\ell,\,\nu}_\tau}(S_\ell(f))\\
&&\quad\quad+\sum_{k'=\ell +1}^\fz\sum_{\tau\in I_{k'}}\sum_{\nu=1}^{N(k',\,\tau)}
\mu(Q_\tau^{k',\,\nu})[S_k(I-S_\ell)\wz D_{k'}](x,\,y^{k',\,\nu}_\tau)
 D_{k'}(f)(y^{k',\,\nu}_\tau).
\end{eqnarray*}
For any $ y_\tau^{\ell,\,\nu}\in Q_\tau^{\ell,\,\nu} $, since
$|S_\ell(f)(z)|\ls G_\ell(f)(y_\tau^{\ell,\,\nu})$
for any $z \in Q_\tau^{\ell,\,\nu}$, we further obtain
$$
|m_{Q^{\ell,\,\nu}_\tau}(S_\ell(f))|\ls G_\ell(f)(y_\tau^{\ell,\,\nu}).$$
Moreover, for any $k'\ge \ell\ge k$ and $y_\tau^{k',\,\nu}\in
Q_\tau^{k',\,\nu}$, by an argument similar to that used in Case 1 of the
proof of \cite[Proposition 5.11]{hmy1}, we have
\begin{eqnarray*}
&&|[S_k(I-S_\ell)\wz D_{k'}](x,\,y^{k',\,\nu}_\tau)|
\ls2^{-(k'-k)\ez}\frac1{V_{2^{-k}}(x)+V(x,\,y^{k',\,\nu}_\tau)}
\lf(\frac{2^{-k}}{2^{-k}+d(x,\,y^{k',\,\nu}_\tau)}\r)^{\ez}.
\end{eqnarray*}
We omit the details here. Thus, for any $ r\in(n/(n+\ez),p)$, we have
\begin{eqnarray*}
&&|S_k(I-S_\ell)f(x)|\\
&&\quad\ls\sum_{k'=\ell}^\fz\sum_{\tau\in
I_{k'}}\sum_{\nu=1}^{N(k',\,\tau)}
\mu(Q_\tau^{k',\,\nu})2^{-(k'-k)\ez}
\frac1{V_{2^{-k}}(x)+V(x,\,y^{k',\,\nu}_\tau)}\\
&&\quad\quad\times
\lf(\frac{2^{-k}}{2^{-k}+d(x,\,y^{k',\,\nu}_\tau)}\r)^{\ez}
 \inf_{z\in
Q_\tau^{k',\,\nu}}G_\ell(f)(z)\\
&&\quad\ls\sum_{k'=\ell}^\fz2^{-(k'-k)(\ez+n-n/r)}
\lf\{\hl\lf(\sum_{\tau\in I_{k'}}\sum_{\nu=1}^{N(k',\,\tau)}
[G_{\ell}(f)]^r\chi_{Q_\tau^{k',\,\nu}}\r)(x)\r\}^{1/r}\\
&&\quad\ls\sum_{k'=\ell}^\fz2^{-(k'-\ell)(\ez+n-n/r)} \lf\{\hl\lf(
[G_\ell(f)]^r \r)(x)\r\}^{1/r}
\ls \lf\{\hl\lf( [G_\ell(f)]^r \r)(x)\r\}^{1/r}.
\end{eqnarray*}
This implies that
$$I_3\ls\lf\|\lf\{\hl\lf( [G_\ell(f)]^r\r)\r\}^{1/r}\r\|_{L^p(\cx)}
\ls\|G_\ell(f)\|_{L^p(\cx)},
$$
which completes the proof of Theorem \ref{t3.3}.
\end{pffive}

Now we turn to the proof of Theorem \ref{t3.4}.

\newtheorem{pfsix}{\it Proof of Theorem \ref{t3.4}}
\renewcommand\thepfsix{}

\begin{pfsix}\rm
Let $\ell\in\zz$.
Then it is quite standard (see, for example,
the proof of \cite[Theorem 2.21]{hmy1}) to prove that for any
$f\in H^{p,\,\fz}_{\ell}(\cx)$,
$$\|G_\ell(f)\|_{L^p(\cx)}
\ls \|f\|_{H^{p,\,\fz}_\ell(\cx)}$$
with a constant independent of $\ell$ and $f$.
Observe that if $f \in(\lip_\ell(1/p-1,\cx))'$,
then $f\in (\cg^\ez_0(\bz,\,\gz))'$ for $\ez\in(n(1/p-1),1]$ and
$\bz,\gz\in(n(1/p-1),\,\ez)$.
Then by Corollary \ref{c3.1}, we have $H^{p,\,\fz}_\ell(\cx)\subset H^p_\ell(\cx)$.
We omit the details.

On the other hand, if $f \in (\cg^\ez_0(\bz,\,\gz))'$ and
$\|G_\ell(f)\|_{L^p(\cx)}<\fz$, then by Theorem \ref{t3.3}, we have
$f-S_\ell f\in H^{p,\,\fz}(\cx)$ and
$$\|f-S_\ell f\|_{H^{p,\,\fz}(\cx)}\sim
\|S^+(f-S_\ell)\|_{L^p(\cx)}
\ls\|G_\ell(f)\|_{L^p(\cx)}.$$
Write
\begin{equation*}
S_\ell(f)=\sum_{\tau\in I_\ell}\sum_{\nu=1}^{N(\ell,\,\nu)}
S_\ell(f)\chi_{Q^{\ell,\,\nu}_\tau}
=\sum_{\tau\in I_\ell}\sum_{\nu=1}^{N(\ell,\,\nu)}
\lz^{\ell,\,\nu}_\tau a^{\ell,\,\nu}_\tau,
\end{equation*}
where $\lz^{\ell,\,\nu}_\tau\equiv[\mu(Q^{\ell,\,\nu}_\tau)]^{1/p}
\|S_\ell(f)\|_{L^\fz(Q^{\ell,\,\nu}_\tau)}$ and
$a^{\ell,\,\nu}_\tau\equiv(\lz^{\ell,\,\nu}_\tau)^{-1}
 S_\ell(f)\chi_{Q^{\ell,\,\nu}_\tau}.$
It is easy to see that
$\{a^{\ell,\,\nu}_\tau\}_{\tau\in I_\ell,\,\nu=1,\,\cdots,\,N(\ell,\,\nu)}$
are
$(p,\,\fz)_{\ell}$-atoms.
Since for any $z\in Q_\tau^{\ell,\,\nu}$,
$$\|S_\ell(f)\|_{L^\fz(Q^{\ell,\,\nu}_\tau)}\ls G_\ell(f)(z),$$
  we have
$$\lf(\lz^{\ell,\,\nu}_\tau\r)^p \ls\int_{Q^{\ell,\,\nu}_\tau}
|G_{\ell}(f)(z)|^p\,d\mu(z),$$
which implies that
$$\sum_{\tau\in I_{\ell}}\sum_{\nu=1}^{N(\ell,\,\nu)}\lf(\lz^{\ell,\,\nu}_\tau\r)^p
\ls\sum_{\tau\in I_{\ell}}\sum_{\nu=1}^{N(\ell,\,\nu)}\int_{Q^{\ell,\,\nu}_\tau}
|G_\ell(f)(z)|^p\,d\mu(z)
\ls\|G_\ell(f)\|_{L^p(\cx)}^p.$$
Thus $S^+_\ell f\in H^p_\ell(\cx)$ and
$\|S^+_\ell f\|\ls\|G_\ell(f)\|_{L^p(\cx)}$,
which implies that
$f\in H^{p,\,\fz}_\ell(\cx)$ and
$\|f\|_{H^{p,\,\fz}_\ell(\cx)}\ls\|G_\ell(f)\|_{L^p(\cx)}$.
This finishes the proof of Theorem \ref{t3.4}.
\end{pfsix}

{\bf Acknowledgements} Dachun Yang would like to thank Professor
Pascal Auscher for some useful discussions on the subject of this paper.
Part of this paper was written during the stay of
Yuan Zhou in University of Jyv\"askyl\"a and he
would like to thank Professor  Pekka Koskela
and Professor Xiao Zhong for their support and hospitality.
Both authors would also like to thank the referee for his/her many valuable remarks
which improve the presentation of the paper.

\end{document}